\documentclass[12pt]{amsart}
\usepackage{latexsym,amsmath,amsfonts,amscd,amssymb}
\newtheorem{thm}{Theorem}[section]
\newtheorem{lemma}[thm]{Lemma}
\newtheorem{prop}[thm]{Proposition}
\newtheorem{cor}[thm]{Corollary}

\theoremstyle{definition}

\theoremstyle{remark}
\newtheorem{remark}[thm]{Remark}

\numberwithin{equation}{section}

\newcommand{\bV}{\bigwedge V}
\newcommand{\bw}{\wedge}
\newcommand{\bbw}{\bigwedge}
\newcommand{\la}{\langle}
\newcommand{\ra}{\rangle}
\newcommand{\too}{\longrightarrow}

\newcommand{\surj}{\twoheadrightarrow}
\newcommand{\inc}{\hookrightarrow}
\newcommand{\bd}{\partial}
\newcommand{\x}{\times}
\newcommand{\ox}{\otimes}

\newcommand{\isom}{\stackrel{\simeq}{\too}}

\newcommand{\coker}{\text{coker}}

\newcommand{\End}{\text{End}}

\newcommand{\Aut}{\text{Aut}}
\newcommand{\Hom}{\text{Hom}}
\newcommand{\Sym}{\text{Sym}}

\newcommand{\Sp}{\text{Sp}}
\newcommand{\SU}{\text{SU}}
\newcommand{\Id}{\text{Id}}

\newcommand{\im}{\text{im}}

\newcommand{\cC}{{\mathcal C}}

\newcommand{\cM}{{\mathcal M}}

\newcommand{\cN}{{\mathcal N}}

\newcommand{\cR}{{\mathcal R}}

\newcommand{\cU}{{\mathcal U}}

\newcommand{\cW}{{\mathcal W}}

\renewcommand{\AA}{{\mathbb A}}
\newcommand{\CC}{{\mathbb C}}

\newcommand{\NN}{{\mathbb N}}
\newcommand{\PP}{{\mathbb P}}
\newcommand{\QQ}{{\mathbb Q}}

\newcommand{\ZZ}{{\mathbb Z}}

\renewcommand{\a}{\alpha}
\renewcommand{\b}{\beta}

\newcommand{\g}{\gamma}

\newcommand{\q}{\psi}

\newcommand{\G}{\Gamma}

\begin{document}

\title[Rational homotopy type of a moduli space]{On the rational
homotopy type of a moduli space of vector bundles over a curve}

\subjclass[2000]{Primary: 14H60; Secondary: 55P62.}
\keywords{Moduli space, rational homotopy, symplectic group.}

\author[I. Biswas]{Indranil Biswas}
\address{School of Mathematics, Tata Institute of Fundamental
Research, Homi Bhabha Road, Bombay 400005, India}
\email{indranil@math.tifr.res.in}

\author[V. Mu{\~n}oz]{Vicente Mu{\~n}oz}

\address{Departamento de Matem\'aticas,
CSIC, Serrano 113 bis, 28006 Madrid, Spain}

\address{Facultad de Matem{\'a}ticas, Universidad Complutense
de Madrid, Plaza de Ciencias 3, 28040 Madrid, Spain}

\email{vicente.munoz@imaff.cfmac.csic.es}

\date{}

\begin{abstract}
We study the rational homotopy of the moduli space $\cN_X$ of
stable vector bundles of rank two and fixed determinant of odd
degree over a compact connected Riemann surface $X$ of genus $g$,
with
$g\geq 2$. The symplectic group $\text{Aut}(H_1(X, {\mathbb Z}))
\,\cong\, \Sp(2g, {\mathbb Z})$ has a natural action on the
rational homotopy groups $\pi_n(\cN_X){\ox}_{\mathbb Z}\QQ$. We
prove that this action extends to an action of $\Sp(2g,\CC)$ on
$\pi_n(\cN_X){\ox}_{\mathbb Z} \CC$. We also show that
$\pi_n(\cN_X){\ox}_{\mathbb Z} \CC$ is a non--trivial
representation of $\Sp(2g,\CC)\,\cong\, \text{Aut}
(H_1(X, {\mathbb C}))$ for all $n\geq 2g-1$. In
particular, $\cN_X$ is a rationally hyperbolic space. In the
special case where $g=2$, for each $n\, \in\, {\mathbb N}$,
we compute the leading $\Sp(2g,\CC)$--representation occurring
in $\pi_n(\cN_X)\ox_\ZZ \CC$.
\end{abstract}

\maketitle

\section{Introduction} \label{sec:introduction}

Moduli spaces of vector bundles over curves have been studied from
various points of view. The aim here is to initiate investigations
of their rational homotopy groups.

Let $X$ be an irreducible smooth projective curve, defined over
$\mathbb C$, of genus $g$, with $g\geq 2$. Fix a holomorphic line
bundle $L_0$ over $X$ of degree one, and consider the moduli space
$\cN_X$ of stable vector bundles $E\to X$ of rank two with
$\bigwedge^2 E\cong L_0$. This moduli space $\cN_X$ is an irreducible
smooth complex projective variety of complex dimension $3g-3$ (see
\cite{Ne2}).

The mapping class group of $X$ acts in a natural way on the
cohomology algebra $H^*(\cN_X, \, {\QQ})$
of $(\cN_X$. This action
actually factors through an action of the symplectic group
$\text{Aut}(H_1(X, {\mathbb Z}))\, \cong\, \Sp(2g,\ZZ)$,
which is a quotient of the mapping class group. Moreover,
the descended action of $\Sp(2g,\ZZ)$ on $H^*(\cN_X,\QQ)$
extends to an action of
$\text{Aut}(H_1(X, {\mathbb C}))\,\cong\,\Sp(2g,\CC)$ on
$H^*(\cN_X,\CC)$. On the other hand, using the fact that
$\cN_X$ is simply connected,
the mapping class group acts naturally on the homotopy groups
$\pi_*(\cN_X)$. Therefore, the mapping class group acts
on $\pi_*(\cN_X){\ox}_{\mathbb Z} \QQ$.

Fix a symplectic basis of $H_1(X, {\mathbb Z})$. Using this
basis $\text{Aut}(H_1(X, {\mathbb Z}))$ (respectively,
$\text{Aut}(H_1(X, {\mathbb C}))$) gets identified with
$\Sp(2g,\ZZ)$ (respectively, $\Sp(2g,\CC)$).

Our first main result is the following (see Theorem \ref{thm:main2}):

\begin{thm} \label{thm:main2-intro}
The action of the
mapping class group on the rational homotopy groups
$\pi_*(\cN_X){\ox}_{\mathbb Z} \QQ$ factors through an action of the
symplectic group ${\rm Sp}(2g,\ZZ)$. This descended action
of ${\rm Sp}(2g,\ZZ)$ extends to
an action of ${\rm Sp}(2g,\CC)$ on $\pi_*(\cN_X)\ox_{\mathbb Z} \CC$.
\end{thm}

We shall prove this theorem in Section \ref{sec:G-equivariant}
using the formality of $\cN_X$
and endowing the minimal model of $\cN_X$ with an
action of $\Sp(2g,\CC)$.

In Sections \ref{sec:moduli} and \ref{sec:general}, we study
the $\Sp(2g,\CC)$--representations
$\pi_n(\cN_X){\ox}_{\mathbb Z} \CC$. In the special
case of $g=2$, we compute the leading representation
for each $n\geq 2$ (Theorem \ref{thm:leading}).

In the general case where $g\,\geq\,2$, we find some non--trivial
irreducible $\Sp(2g,\CC)$--representations
contained in $\pi_n(\cN_X)\ox_\ZZ \CC$
for each $n\geq 2g$ (see Theorem \ref{thm:leading} for the case of
$g=2$ and Theorem \ref{thm:g>2} for $g>2$). We have the following
result.

\begin{thm} \label{thm:main}
Take any integer $n$ with $n\geq 2g$. The ${\rm
Sp}(2g,\CC)$--module $\pi_n(\cN_X)\ox_\ZZ \CC$ is non--trivial. So
the action of ${\rm Sp}(2g,\ZZ)$ on the rational homotopy groups
$\pi_n(\cN_X)\ox_\ZZ \QQ$ is non--trivial, and the action of the
mapping class group on $\pi_n(\cN_X)$ is non--trivial.
\end{thm}

A connected simply connected finite CW--complex $Z$ (e.g.\ a
compact $1$--connected manifold) is said to be \textit{rationally
elliptic} if the total dimension of the rational homotopy groups
is finite, or in other words,
 \[
 \sum_{n\in {\mathbb N}} \dim \pi_n(Z)\ox_\ZZ\QQ\,<\, \infty\, .
 \]
Otherwise, $Z$ is called \textit{rationally hyperbolic} (see
\cite{Felix}). If $Z$ is rationally elliptic of dimension $N$,
then $\pi_n(Z)\ox_\ZZ\QQ=0$ for all $n\geq 2N$ (equivalently,
$\pi_n(Z)$ are torsion for $n\geq 2N$). On the other hand, if $Z$
is rationally hyperbolic of dimension $N$, then
 \[
 f(k)\,=\, \sum_{i=1}^{N-1} \dim \pi_{k+i}(Z){\ox}_{\mathbb Z} \QQ
 \]
grows faster than any polynomial in $k$. This dichotomy is
discussed in \cite{Felix}.

A byproduct of Theorem \ref{thm:main} is the following corollary.

\begin{cor}\label{cor:hyperbolic}
The moduli space
$\cN_X$ is rationally hyperbolic for all $g\, \geq\, 2$.
\hfill $\Box$
\end{cor}

\section{The moduli space $\cN_X$} \label{sec:introduction2}

Let $X$ be an irreducible smooth complex projective curve of genus
$g\geq 2$. Fix a holomorphic line bundle $L_0$ over $X$ of degree
one. Let $\cN_X$ denote the moduli space of stable vector bundles
over $X$ of rank two and $\det(E)\, =\, \bigwedge^2 E\, = \, L_0$.
This moduli space $\cN_X$ is an irreducible smooth complex
projective variety of complex dimension $3g-3$ (see \cite{Ne2}).
In particular, it is a compact connected $C^\infty$ (real)
manifold of dimension $6g-6$. The complex structure of $X$ endows
$\cN_X$ with a natural K{\"a}hler structure \cite{AB}.

If we take any holomorphic line bundle $L_1$ over $X$ of odd
degree $2d+1$, then there is a holomorphic line bundle $\mu$ over
$X$ of degree $d$ such that $L_1\cong L_0\ox \mu^2$. Therefore,
the map defined by
$E\longmapsto E\ox \mu$ is an algebraic isomorphism from
$\cN_X$ to the moduli space of stable vector bundles of
rank two over $X$ with determinant $L_1$. In particular, the
isomorphism class of the variety $\cN_X$ is independent of the
choice of the line bundle $L_0$.

The diffeomorphism class of the real manifold $\cN_X$ is
independent of the complex structure of $X$. This can be seen as
follows. Fix a point $x_0\in X$, and set $X'=X\setminus \{x_0\}$
to be the complement. Choosing a point $x'\in X'$, consider
the subset
$$
 \Hom^0(\pi_1(X',x')\, , \SU(2)) \,
 \subset\, \Hom (\pi_1(X',x')\, , \SU(2))
$$
parametrizing all homomorphisms from the fundamental group
$\pi_1(X',x')$ to ${\rm SU}(2)$ satisfying the condition that the image
of the conjugacy class in $\pi_1(X',x')$ corresponding to the free
homotopy class of oriented loops in $X'$ around $x_0$ (with
anticlockwise orientation) is $-\Id$. Let
 \begin{equation}\label{eqn:Rg}
 {\cR}_g\, :=\, \Hom^0(\pi_1(X',x')\, , \SU(2))/\SU(2)
 \end{equation}
be the quotient space for the adjoint action of $\SU(2)$
on itself. It is easy to see that
$\cR_g$ is a connected compact $C^\infty$ manifold of dimension
$6g-6$ (see \cite{Ne0}). Given any homomorphism $\rho\, \in\,
\Hom^0(\pi_1(X',x'),{\rm SU}(2))$, the corresponding flat vector
bundle over the Riemann surface $X'$ extends to $X$ as a
holomorphic vector bundle with a logarithmic connection which has
residue $-\frac{1}{2}\Id$ at $x_0$ (see \cite{De}). The underlying
holomorphic vector bundle $E_\rho$ on $X$ is stable and
$\det(E_\rho) \,=\, {\mathcal O}_X(x_0)$. Sending any $\rho$ to
$E_\rho$ we obtain a diffeomorphism of ${\cR}_g$ with $\cN_X$ for
$L_0\, =\, {\mathcal O}_X(x_0)$ (cf.\ \cite{NS}). Therefore, the
diffeomorphism class of the real manifold $\cN_X$ is independent
of the complex structure of $X$ (it depends only on the genus of
$X$).

It is easy to see that the manifold $\cR_g$ in eqn. \eqref{eqn:Rg} is
simply connected \cite[Corollary 2]{Ne0}. By \cite[Theorem
9.10]{AB}, the cohomology ring $H^*(\cN_X,\ZZ)$ is torsionfree,
and by \cite[Proposition 9.13]{AB}, we have $H^2(\cN_X,\,
\ZZ)=\ZZ$. Consequently, the variety $\cN_X$ has a natural
polarization. Henceforth, we shall denote by $\alpha$ the natural
positive (i.e., ample) generator of $H^2(\cN_X,\, \ZZ)$.

Next we will describe an action of the mapping class group on the
cohomology of $\cN_X$. For that purpose, consider the
moduli space ${\cM}^1_g$ parametrizing all isomorphism classes of
one--pointed compact Riemann surfaces $(Y,y)$ of genus $g$ with
$\Aut(Y,y) \, =\, e$ (i.e., $Y$ does not have any nontrivial
automorphism that fixes the marked point $y\, \in\, Y$).
This moduli space is a smooth irreducible
quasiprojective variety of dimension $3g-2$ defined over the field
$\CC$. Given any $(Y,y)$, there is a natural choice of a
holomorphic line bundle of degree one over $Y$, namely
${\mathcal O}_Y(y)$.
There is a universal family of Riemann surfaces
 \begin{equation}\label{eqn:uf}
 p\, \colon\, {\cC}_g \, \too\, {\cM}^1_g
 \end{equation}
and a holomorphic section $h\, \colon\, {\cM}^1_g \,
\too\,{\cC}_g$ giving the marked point. Let
 \begin{equation}\label{eqn:uh}
 P\, \colon\, \widetilde{\cN}\, \too\, {\cM}^1_g
 \end{equation}
be the family of moduli spaces of stable vector bundles of rank
two with fixed determinant corresponding to the family of Riemann
surfaces in eqn. \eqref{eqn:uf}. For any one--pointed Riemann surface
$\underline{x}=(Y,y) \, \in\, \cM^1_g$, the fiber
$P^{-1}(\underline{x})$ is the moduli space ${\cN}_Y$
parametrizing all stable vector bundles over $Y$ of rank two and
determinant ${\mathcal O}_Y(y)$.

Fix a base point $\underline{x}_0\,=\,(X\, , x_0) \, \in\, \cM^1_g$
of the moduli space.
Let $G_{\ZZ}$ (respectively, $G_{\CC}$) denote the group of all
automorphisms of $H^1(X,\, {\ZZ})$ (respectively, $H^1(X,\,
{\CC})$) preserving the symplectic pairing given by the cup
product. Choosing a symplectic basis of $H^1(X,\, {\ZZ})$, the
groups $G_{\ZZ}$ and $G_{\CC}$ get identified with $\Sp(2g,
{\ZZ})$ and $\Sp(2g, {\CC})$ respectively.

\noindent
\textbf{Convention.}\, In the sequel, we will
interchange $G_{\ZZ}$ (respectively, $G_{\CC}$)
and $\Sp(2g, {\ZZ})$ (respectively, $\Sp(2g, {\CC})$).

Consider the local system $R^{1}p_*\underline{\underline{\mathbb Z}}$
on ${\cM}^1_g$, where $p$ is the projection in eqn. \eqref{eqn:uf},
and $\underline{\underline{\ZZ}}$ is the constant local system on
${\cC}_g$ with stalk $\ZZ$. Using its monodromy, the group
$G_{\ZZ}$ is a quotient of the fundamental group
 $$
 \G^1_g \, :=\, \pi_1({\mathcal M}^1_g,\, \underline{x}_0)\, .
 $$
This group $\G^1_g$ is known as the \textit{mapping class group},
and the kernel of the projection of $\G^1_g$ to $G_{\ZZ}$ is known
as the \textit{Torelli group}.

Actually, the mapping class group has a natural
action on the moduli space $\cN_X\,=\,
P^{-1}(\underline{x}_0)$. To see this action, note
that using the earlier mentioned
identification $\cN_X\,=\, \cR_g$ (defined in eqn. \eqref{eqn:Rg}),
the fiber
bundle $P$ in eqn. \eqref{eqn:uh} has a natural flat connection (this
flat connection is not holomorphic). The monodromy of this
flat connection gives an action of
$\G^1_g \,=\,\pi_1({\cM}^1_g,\, \underline{x}_0)$
on $\cN_X$; more details can be found in \cite{Bi}.

The action of $\G^1_g$ on $H^i(\cN_X,\, {\mathbb Z})$ induced by
the above action of $\G^1_g$ on $\cN_X$ evidently coincides with
the monodromy representation of the local system $R^iP_*
\underline{\mathbb Z}$ on ${\cM}^1_g$, where $\underline{\mathbb
Z}$ is the constant local system on $\widetilde{\cN}$ with stalk
$\mathbb Z$.

\begin{prop}\label{prop:action}
The action of the mapping class group on the cohomology algebra
$H^*(\cN_X, \, {\QQ})$ factors through an action of the symplectic
group $G_{\mathbb Z}\, =\,
{\rm Sp}(2g,\ZZ)$. Moreover, this action of ${\rm Sp}(2g,\ZZ)$
on $H^i(\cN_X,\QQ)$ extends to an action of $G_{\mathbb C}\, =\,{\rm
Sp}(2g,\CC)$ on $H^i(\cN_X,\CC)$.
\end{prop}

\begin{proof}
The cohomology algebra
$H^*(\cN_X,\, {\QQ})$ is generated by the K\"unneth components of
the second Chern class of the adjoint bundle of
a universal vector bundle over $X\times {\cN}_X$ (see \cite{Ne3, AB}
and also Section \ref{sec:cohomology}). Note that although
there is no unique universal bundle
over $X\times {\cN}_X$, any two universal bundles differ
by tensoring with a line bundle pulled back from ${\cN}_X$. Therefore,
the universal adjoint bundle is unique. Consequently,
the local system $\bigoplus_{i\geq 0} R^{i}P_*\underline{{\CC}}$
on ${\cM}^1_g$,
where $\underline{\CC}$ is the constant local system on
$\widetilde{\cN}$ with stalk $\CC$, is a quotient of some local
system on ${\cM}^1_g$ of the form
 $$
 {\cW}\, :=\, \bigoplus_{j=1}^\ell\big(\big( \bigoplus_{i=0}^2
 R^{i}p_*\underline{\underline{\mathbb C}} \big)^{\oplus
 a_j}\big)^{\ox b_j}\, ,
 $$
where $\ell ,a_j,b_j\, \in\, {\NN}$, the map $p$ is the projection in
eqn. \eqref{eqn:uf} and $\underline{\underline{\CC}}$ is the constant
local system on ${\cC}_g$ with stalk $\CC$. In other words, we have
a surjective homomorphism of local systems
 \begin{equation}\label{eqn:surj.}
 {\cW}\, \longrightarrow\, \bigoplus_{i\geq 0} R^{i}P_*\underline{{\CC}}
 \, \longrightarrow\, 0\, .
 \end{equation}
Both $R^{0}p_*\underline{\underline{\mathbb C}}$ and
$R^{2}p_*\underline{\underline{\mathbb C}}$ are constant local
systems on ${\cM}^1_g$, and the monodromy of the local system
$R^{1}p_*\underline{\underline{\mathbb C}}$\ , by definition,
factors through $G_{\mathbb Z}$. Consequently, the monodromy
representation
 \begin{equation}\label{eqn:mo.}
 \G^1_g\, \too\, \text{Aut}({\cW}_{\underline{x}_0})
 \end{equation}
of the mapping class group for the local system $\cW$
on ${\cM}^1_g$ factors through $G_{\mathbb Z}$.
Hence, the Torelli group is in the kernel of the monodromy
representation
 \begin{equation}\label{eqn:mh}
 \G^1_g\, \too\, \prod_{i\geq 0} \text{Aut}(
 (R^{i}P_*\underline{{\mathbb C}})_{\underline{x}_0})
 \end{equation}
of the mapping class group for the quotient local system in
eqn. \eqref{eqn:surj.}. Therefore, the homomorphism in eqn.
\eqref{eqn:mh}
factors through the quotient $G_{\mathbb Z}$ of $\G^1_g$.

To prove that the action of $\Sp(2g,\ZZ)$ on $H^i(\cN_X,\QQ)$
extends to an action of $G_{\mathbb C}\, =\,
\Sp(2g,\CC)$ on $H^i(\cN_X,\CC)$, first
note that the monodromy representation
$$
G_{\mathbb Z}\, \too\,
\text{Aut}({\cW}_{\underline{x}_0})
$$
in eqn. \eqref{eqn:mo.} extends to
a homomorphism from $\Sp(2g,\CC)$. The kernel of the surjective
homomorphism
 \begin{equation}\label{eqn:ker.}
 {\cW}_{\underline{x}_0}\, \longrightarrow\, \bigoplus_{i\geq 0}
 (R^{i}P_*\underline{{\CC}})_{\underline{x}_0}
 \end{equation}
obtained from eqn. \eqref{eqn:surj.} is preserved by $G_{\mathbb Z}$.
On the other hand, $\Sp(2g,\ZZ)$ is Zariski dense in $\Sp(2g,\CC)$
(see \cite{Bo}). Hence the kernel of the homomorphism in
eqn. \eqref{eqn:ker.} is preserved by the action of $\Sp(2g,\CC)$ on
${\cW}_{\underline{x}_0}$. Consequently, the action of
$\Sp(2g,\CC)$ on ${\cW}_{\underline{x}_0}$ induces an action of
$\Sp(2g,\CC)$ on the quotient in eqn. \eqref{eqn:ker.}. This completes
the proof of the proposition.
\end{proof}

\section{Cohomology ring of $\cN_X$} \label{sec:cohomology}

Let us recall the known description of the cohomology ring
$H^*(\cN_X,\QQ)$ of the moduli space $\cN_X$ (see~\cite{King, ST,
Quantum}). Consider a universal bundle $\cU \to X \x \cN_X$. Let
$\End_0 (\cU) \to X \x \cN_X$ be the adjoint vector bundle (we recall
that $\End_0 (\cU)\subset \text{End} (\cU)$ is subbundle of rank
three given by the trace--free endomorphisms of
the fibers of $\cU$). The
K\"unneth decomposition of the second Chern class
$c_2(\End_0(\cU)) \in H^4(X \x \cN_X,\ZZ)$ can be written as
\begin{equation} \label{eqn:qu2}
c_2(\End_0 \, (\cU))=2 [X] \ox \a + 4 \, \q -1 \ox \b\, ,
\end{equation}
where $\b\in H^4(\cN_X,\ZZ)$,
$[X]\in H^2(X,\ZZ)$ denotes the fundamental class of the
Riemann surface $X$,
$\a\in H^2(\cN_X,\ZZ)$ as before is the positive generator of
$H^2(\cN_X,\ZZ)= \ZZ$, and $\q\in
H^1(X,\ZZ)\ox_\ZZ H^3(\cN_X,\ZZ)$. Let $\{ c_1,\ldots, c_{2g}\}$
be a symplectic basis of $H^1(X,\ZZ)$, which means that $c_i \cup
c_{i+g}=[X]$ for all $1\leq i\leq g$, and $c_j \cup c_k=0$ for all
$j,k$ with $|j-k|\neq g$. It is known that $\q\,=\,
\sum_{i=1}^{2g} c_i \ox \g_i$,
where $\{\g_1,\ldots, \g_{2g}\}$ is a basis for $H^3(\cN_X,\ZZ)$;
see \cite{MN}. In other words, $\q$ gives an isomorphism
\begin{equation}\label{eqn:h3}
H^1(X,\ZZ)\, =\, H^1(X,\ZZ)^* \,\too\, H^3(\cN,\ZZ).
 \end{equation}

The elements $\a, \b$ and
$\g_i$, $1\leq i \leq 2g$, together generate
$H^*(\cN_X,\QQ)$ as an algebra \cite{Ne3, AB, Thaddeus}. We can
rephrase this as saying that there exists an epimorphism
\begin{equation}\label{eqn:epim}
F \colon \bbw(\a,\g_1,\ldots, \g_{2g},\b) :=
\QQ[\a,\b ]\ox \bw(\g_1,\ldots, \g_{2g}) \surj H^*(\cN_X,\QQ),
\end{equation}
where $\deg(\a)=2$, $\deg(\b)=4$ and $\deg(\g_i)=3$, $1\leq i \leq
2g$. Here $\bbw$ means the free graded algebra generated by the
given elements, which is the tensor product of the symmetric
algebra on the even degree elements and the exterior algebra on the
odd degree elements.

We shall denote by $W$ the standard $\mathbb Q$--representation of
$G_{\mathbb Z}\, =\, \Sp(2g,\ZZ)$, so
  $$
  H^1(X,\QQ)\,\cong\, W\, .
  $$
We noted in Section \ref{sec:introduction2}
that the monodromy action of $\G_g^1$ on $H^*(\cN_X,\QQ)$
factors through an action of $\Sp(2g,\ZZ)$. It is easy to see that
this action fixes both
$\alpha$ and $\beta$, and furthermore,the isomorphism in
eqn. \eqref{eqn:h3} is $\Sp(2g,\ZZ)$--equivariant. Therefore,
 $$
 H^3(\cN_X,\QQ)\,\cong\, H^1(X,\QQ)^*\,\cong\, W^* \,\cong\, W
 $$
as $\Sp(2g,\ZZ)$--representations.

Let
 $$
 H_I^*(\cN_X,\QQ)\, \subset\, H^*(\cN_X,\QQ)
 $$
be the subalgebra fixed pointwise by the action of $\Sp(2g,\ZZ)$.
The epimorphism in eqn. \eqref{eqn:epim} is
$\Sp(2g,\ZZ)$--equivariant, and $\Sp(2g,\ZZ)$ is Zariski dense in
the reductive group $\Sp(2g,{\mathbb C})$ \cite{Bo}. Using these
we conclude that the invariant part $H_I^*(\cN_X,\QQ)$ is
generated by $\a$, $\b$ and $\g=-2 \sum_{i=1}^g \g_i\g_{i+g}$ (the
factor of $-2$ is for convenience, to be in accordance with the
existing literature). Then the epimorphism $F$ in eqn.
\eqref{eqn:epim} gives an epimorphism
 \begin{equation}
  \QQ[\a,\b,\g ]  \,\surj\, H^*_I(\cN_X,\QQ)\, ,
 \label{eqn:qu4}
  \end{equation}
where $\deg(\a)=2$, $\deg(\b)=4$ and $\deg(\g)=6$. Hence we may
write
 \begin{equation}\label{eqn:invar}
 H_I^*(\cN_X,\QQ)= \frac{\QQ [\a, \b, \g]}{I_g}\, ,
 \end{equation}
where $I_g$ is an ideal of relations satisfied by $\a$, $\b$ and
$\g$.

For each $0 \leq k \leq g$, the primitive component of $\wedge^k
W$ is defined as
 $$
 \bw_0^k W = \ker (\g^{g-k+1} \colon \bw^k W \too \bw^{2g-k+2} W).
 $$
The spaces $\bw^k_0 W$ are irreducible
$\Sp(2g,\ZZ)$--representations.

The descriptions of the ideal
$I_g$ and the cohomology ring $H^*(\cN_X,\QQ)$ are given in the
following proposition.

\begin{prop}[\cite{King, ST}]
\label{prop:cohomology} Define $q^1_0=1$, $q^2_0=0$, $q^3_0=0$ and
then recursively, for all $r \geq 1$,
 $$
 \left\{ \begin{array}{l} q_{r+1}^1 = \a q_r^1 + r^2 q_r^2\, ,
 \\ q_{r+1}^2 = \b q_r^1 + \frac{2r}{r+1}  q_r^3\, ,
 \\ q_{r+1}^3 =  \g  q_r^1\, .
 \end{array} \right.
 $$
Then $I_g=(q^1_g,q^2_g,q^3_g) \subset \QQ[\a,\b,\g]$, for all $g
\geq 1$. Note that $\deg(q_g^1)=2g$, $\deg(q_g^2)=2g+2$ and
$\deg(q_g^3)=2g+4$. Moreover the ${\rm
Sp}(2g,{\ZZ})$--decomposition of $H^*(\cN_X,\QQ)$ is
 \begin{equation} \label{eqn:decomposition}
 H^*(\cN_X,\QQ)= \bigoplus_{k=0}^{g-1} \bw_0^k W
 \ox \frac{\QQ [\a, \b, \g]}{I_{g-k}}\, . \quad \Box
 \end{equation}
\end{prop}

\begin{lemma} \label{lem:kerF}
The vector space
  $$
  E = \la q_g^1\ra \oplus\la q_g^2\ra \oplus (q_{g-1}^1 \cdot W) \oplus
  (q_{g-2}^1\cdot \bw^2_0 W) \oplus \cdots \oplus  (q_{1}^1\cdot
  \bw^{g-1}_0 W) \oplus
  \bw^g_0 W\, ,
  $$
realized as a subspace
of $\AA\,:= \,\bbw (\alpha, \g_1,\ldots,\g_{2g},\beta)
\,=\, {\mathbb Q}[\alpha\, ,\beta]\otimes \wedge
(\g_1,\ldots,\g_{2g})$ using the identification
$W\, =\, \la\g_1,\ldots,\g_{2g}\ra$, generates
the ideal ${\rm kernel}(F)$ of the map $F$ in eqn. \eqref{eqn:epim}.
\end{lemma}

\begin{proof}
Clearly we have $E\,\subset\, {\rm kernel}(F)$. We will prove the
reverse inclusion
 $$
 {\rm kernel}(F)\,\subset\, I(E)\, ,
 $$
where $I(E)$ is the ideal generated by $E$ in $\AA$.

By Proposition \ref{prop:cohomology}, ${\rm kernel}(F)$ is
generated by $q_{g-k}^i \cdot \bw^k_0 W$, where $i\in [1,3]$ and
$k\in [0,g]$. Note that since $q_{0}^2=0$ and $q_0^3=0$, it suffices
to prove the following two:
\begin{enumerate}
\item $q^2_{g-k} \cdot \bw^k_0  W \, \subset\, I(E)$
for $1\leq k \leq g-1$, and

\item $q^3_{g-k} \cdot \bw^k_0 W \,
\subset\, I(E)$ for $0\leq k \leq g-1$.
\end{enumerate}

We shall use the following inclusions:
 \begin{eqnarray}
 && \gamma \cdot \bw^{j}_0 W \subset I(\bw^{j+1}_0 W), \qquad 0\leq 
j\leq g-1 \, , \label{eqn:extra1} \\
 && \bw^{j+1}_0 W \subset I(\bw^{j}_0 W), \qquad 0\leq j \leq 
g-1\label{eqn:extra2} \,.
 \end{eqnarray}

For proving
eqn. (\ref{eqn:extra1}), first note that $\bw^{j}_0 W$ is an
irreducible ${\rm Sp}(2g,{\ZZ})$--representation, so it is enough
to see that there is a non-zero element in $\gamma \cdot \bw^{j}_0
W$ which lies in $I(\bw^{j+1}_0 W)$. Consider $\g_1\cdots \g_j\in
\bw^j_0 W$. Then
$$
\g\cdot \g_1\cdots \g_j\,=\, - 2 \sum_{i=j+1}^g
\g_1\cdots \g_j\g_{j+1}\g_{j+1+g}\, ,
$$
and $\g_1\cdots \g_{j+1}\in
\bw^{j+1}_0 W$. Therefore $\g\cdot \g_1\cdots \g_j \in
I(\bw^{j+1}_0 W)$, as required.

To prove eqn. (\ref{eqn:extra2}),
we first note that $\bw^{j+1}_0 W$ is an
irreducible ${\rm Sp}(2g,{\ZZ})$--representation and $\g_1\cdots
\g_{j+1}\in \bw^{j+1}_0 W$. Clearly we have $\g_1\cdots \g_j \in
\bw^j_0 W$. Hence it follows that $\g_1\cdots \g_j \g_{j+1}\,\in\,
I(\bw^{j}_0 W)$. This gives the required inclusion.

Using eqn. (\ref{eqn:extra1}), we have that
 $$
 q_{g-k}^3 \cdot \bw^k_0 W = q_{g-k-1}^1\gamma \cdot \bw^k_0 W \subset
 I(q_{g-k-1}^1 \cdot \bw^{k+1}_0 W )\subset I(E)
 $$
for all $0\leq k\leq g-1$. Also, using eqn. (\ref{eqn:extra2}) we have
  \begin{eqnarray*}
  q_{g-k}^2 \cdot \bw^k_0 W &=& \frac{1}{(g-k)^2}(q_{g-k+1}^1 - \alpha
  q_{g-k}^1) \cdot \bw^k_0 W \\&\subset&
  I(  q_{g-k+1}^1 \cdot \bw^{k-1}_0 W \oplus  q_{g-k}^1\cdot \bw^k_0 W)
  \subset I(E)\, ,
  \end{eqnarray*}
 for all $1\leq k\leq g-1$.
\end{proof}

\begin{remark}
  The subspace $E$ in Lemma \ref{lem:kerF} is minimal in the
  sense that no proper subspace of $E$ generates ${\rm kernel}(F)$.
\end{remark}

\section{Minimal models}\label{sec:formal}

Let us recall some definitions and results about minimal models
(see \cite{DGMS,GM} for more details). Let $(A,d)$ be a {\it
differential graded algebra} (in the sequel, we shall just say a
differential algebra). This means that $A$ is a graded (in
non--negative degrees) commutative algebra over a field $K$, of
characteristic zero, and $d\colon A^n\too A^{n+1}$ is a
differential which satisfies the derivation condition which says
that
 $$
 d(a\cdot b) \,=\, (da)\cdot b +(-1)^{\deg (a)} a\cdot (db)\, ,
 $$
where $\deg(a)$ is the degree of $a$. Throughout this article we
shall assume that $K=\CC$, the field of complex numbers.

Morphisms between differential algebras are required to be degree
preserving algebra maps that commute with the differentials.
Given a differential algebra $(A,d)$, we denote by $H^*(A,d)$ its
cohomology. We say that $A$ is {\em connected\/} if
$H^0(A,d)=\CC$, and {\em $1$--connected\/} if, in addition,
$H^1(A,d)=0$.

A differential algebra $(A,d)$ is said to be {\it minimal\/} if
the following two hold:
\begin{enumerate}
 \item[(i)] $A$ is free as a graded algebra, that is, $A\, =\,
\bigwedge V$, where $V=\oplus_{i>0} V^i$ is a graded vector space,
and
 \item[(ii)] there exists a collection of generators
$\{a_\tau\}_{\tau\in I}$ of the algebra $A$, where $I$ is some
well ordered index set, such that $\deg(a_\mu)\leq \deg(a_\tau)$
if $\mu < \tau$ and each $d a_\tau$ is expressed in terms of
preceding $a_\mu$, $\mu<\tau$.
\end{enumerate}

As before, $\bV$ is the tensor product of the symmetric algebra
on the even degree part of $V$ with the exterior algebra
on the odd degree part of $V$

For notational convenience, we shall use the dot ``$\cdot$'' to
denote the product operation on $\bV$.

For any $n$, define $V^{\leq n}\,:=\, \bigoplus_{i\leq
n} V^i$. So $\bV^{\leq n}\, =\, \bbw (\bigoplus_{i\leq n} V^i)$ is
the subalgebra generated by elements of degrees at most $n$. For
any $m$, let $(\bV)^m$ denote the subspace of $\bV$
spanned by all elements of total degree $m$.
Finally, for $k\geq 1$, let $\bbw^{\geq k} V$ denote the
ideal formed by elements which are products of at least $k$
generators. In other words,
 $$
 \bbw\nolimits^+ V\, :=\, (\bV)^{>0}\, =\, \bigoplus_{m>0}(\bV)^m\, ,
 $$
and
 $$
 \bbw\nolimits^{\geq k}\,=\,
 \overset{k-\text{times}}{\overbrace{\bbw\nolimits^+ V
 \cdots \bbw\nolimits^+ V}} \, .
 $$

Note that the condition (ii) in the definition of minimality
implies that $d\colon V\too \bbw^{\geq 2}V$, and hence $d\colon
\bbw^{i} V\too \bbw^{\geq (i+1)} V$, for all $i\geq 1$. Notations
like $(\bbw^{\geq i} V^{<j})^n$ have natural meaning.

Given a differential algebra $(A,d)$,
we shall say that $({\bV},d)$ is a {\it minimal model} of
$(A,d)$ if $(\bV,d)$ is minimal and there
exists a morphism of differential graded algebras $\rho\colon
{({\bV},d)}\too {(A,d)}$ such that the induced homomorphism
of cohomologies
 $$
 \rho^*\,\colon\, H^*(\bV,d)\too H^*(A,d)
 $$
is an isomorphism. Such a homomorphism $\rho$ is called a
\emph{quasi--isomorphism}. Any $1$--connected differential algebra
$(A,d)$ has a minimal model unique up to an isomorphism \cite{DGMS,
GM}.

A {\it minimal model\/} of a connected differentiable manifold $M$
is a minimal model $(\bigwedge V,d)$ for the de Rham complex
$(\Omega^*( M,\CC),d)$ of complex
$C^\infty$ differential forms on $M$. If $M$ is simply
connected, then the dual of the complex homotopy vector
space $\pi_i(M)\ox_\ZZ \CC$ is isomorphic to $V^i$ for any $i>0$
(see \cite{GM}).

A minimal model $({\bV},d)$ is said to be {\it formal} if there is
a morphism of differential algebras
$$
\psi\,\colon\, {({\bV},d)}\,\too\,
(H^*(\bV,d),0)
$$
which induces the identity map on cohomology. This
means that $(\bV,d)$ is the minimal model of the algebra
$(H^*(\bV,d),0)$ with zero differential.

We shall say that a connected differentiable manifold $M$ is
\emph{formal} if its minimal model is formal, or equivalently, the
two differential algebras $(\Omega^*( M,\CC),d)$ and $(H^*(M,
{\mathbb C})\, , 0)$ have the same minimal model. Therefore, if
$M$ is formal and simply connected, then the complex homotopy
groups $\pi_i(M)\ox_\ZZ\CC$ are obtained by computing the minimal
model of $(H^*(M, {\CC})\, ,0)$.

The main result of \cite{DGMS} gives the following strong
topological restriction on the rational homotopy type of K{\"a}hler
manifolds.

\begin{thm}[\cite{DGMS}]\label{prop:DGMS}
Let $M$ be a compact connected K\"ahler manifold.
Then $M$ is formal. \hfill $\Box$
\end{thm}

Therefore, the minimal model of a compact connected K\"ahler
manifold $M$ can be obtained from the minimal model of its
cohomology algebra $(H^*(M, {\CC})\, ,0)$. Moreover, if the
K\"ahler manifold $M$ is simply connected, this process will also
gives us the complex homotopy groups $\pi_i(M)\ox_\ZZ {\CC}$ of
$M$.

We will briefly review a construction of the minimal model of a
differential
algebra $(A,d)$. For simplicity, we shall assume that $(A,d)$ is
$1$--connected. We need to find a graded vector space $V=\oplus_{n
\geq 1} V^n$, a differential $d$, with
$$
d|_{V^n}\colon  V^n \too
\bbw^{\geq 2} V^{\leq (n-1)}\, ,
$$
and a graded linear map
 $$
 \rho\, =\, \sum
 \rho_n\colon  V\,=\, \oplus V^n \,\too\, A\,=\, \oplus A^n
 $$
such that the induced homomorphism $\rho\colon  \bV\too A$
respects the differentials, which means that $\rho\circ d=d\circ
\rho$, and furthermore, the map on cohomology
$$
\rho^*\,\colon\, H^*(\bV,d) \,\too\, H^*(A,d)
$$
is an isomorphism.

We shall construct $V^n$, $\rho_n$ and $d|_{V^n}$, where $n\geq 1$,
using induction on $n$. They will satisfy the following conditions:
 \begin{enumerate}
 \item[(i)] $\rho_n\colon V^n \too A^n$\ ;
 \item[(ii)] $d_n=d|_{V^n}\colon V^n\too \bbw^{\geq 2} V^{\leq (n-1)}$\ ;
 \item[(iii)] $\rho_{\leq (n-1)} \circ d_n= d \circ \rho_n$ on $V^n$,
  where $\rho_{\leq (n-1)} \colon  \bV^{\leq (n-1)} \too A$ is
 induced by the map $\rho_i$, $i\leq n-1$\ ;
 \item[(iv)] $\rho_{\leq n}^*\colon
H^i(\bV^{\leq n},d)\isom H^i(A,d)$ is an isomorphism for $i\leq n$\ ;
 \item[(v)] $\rho_{\leq n}^*\colon  H^{n+1}(\bV^{\leq n},d) \inc
 H^{n+1}(A,d)$ is an injection.
 \end{enumerate}
{}From these conditions it follows that the map $\rho\colon
(\bV,d)\too (A,d)$, constructed using $\rho_n$ on each subspace
$V^n$, is a quasi--isomorphism. Given any $i$, we evidently have
$(\bV)^k=(\bV^{\leq (i+1)})^k$ for all $k\leq i+1$. So $H^i(\bV,d)
\cong H^i(\bV^{\leq (i+1)},d)$. The composition
 $$
 (\bV^{\leq (i+1)},d) \,\inc\, (\bV,d)\, \stackrel{\rho}{\too}\, (A,d)
 $$
equals $\rho_{\leq (i+1)}$. Hence
 $$
 \rho^*\,=\, \rho_{\leq (i+1)}^*\,\colon\,
 H^i(\bV,d) \,\cong\, H^i(\bV^{\leq (i+1)},d) \,\isom\, H^i(A,d)
 $$
is an isomorphism. This proves that $(\bV,d)$ is a minimal model for
$(A,d)$.

To construct $V^n$, $\rho_n$ and $d|_{V^n}$, we start with
$V^1=0$. All conditions (i)---(v) hold trivially, since
$H^1(A,d)=0$.

Now assume that conditions (i)---(v) are satisfied for
all $j\,\in\, [1\, ,n-1]$ with $n-1 \geq
1$; let us see that we can find $V^n$, $\rho_n$ and $d_n$ also
fulfilling these conditions. Take
 $$
 \begin{array}{l}
  V^n= C^n\oplus N^n, \\
  C^n = \coker \left(\rho^*_{\leq (n-1)} \colon  H^n(\bV^{\leq (n-1)},d)
  \inc H^n(A,d) \right), \\
  N^n= \ker \left( \rho_{\leq (n-1)}^*\colon H^{n+1}(\bV^{\leq (n-1)},d)
  \too H^{n+1}(A,d)\right).
 \end{array}
 $$

Define $\rho_n\colon V^n \too A^n$ as follows. First, we introduce
the notation
 \begin{eqnarray*}
 Z^n(A,d) &=& \ker (d \colon A^n\too A^{n+1}), \\
 B^n(A,d) &=&\im (d \colon A^{n-1}\too A^n),
 \end{eqnarray*}
for the spaces of cocycles and coboundaries, respectively. Let
$\imath_1\colon C^n\too H^n(A,d)$ be a linear map which is a
splitting of the projection $H^n(A,d) \surj C^n$. Also, let
$\imath_2\,\colon\, H^n(A,d) \,\too\, Z^n(A,d)$ be a splitting of
the projection $Z^n(A,d) \surj H^n(A,d)$. Let $\imath_3\colon
Z^n(A,d)\inc A^n$ be the inclusion map. Then define
 $$
 \rho_n|_{C^n}\,=\,\imath_3 \circ\imath_2\circ \imath_1\, .
 $$
To define $\rho_n|_{N^n}$, let $\jmath_1\colon N^n \inc
H^{n+1}(\bV^{\leq (n-1)},d)$ be the inclusion. Take a splitting
$\jmath_2\colon  H^{n+1}(\bV^{\leq (n-1)},d)\too Z^{n+1}(\bV^{\leq
(n-1)},d)$ of the obvious projection. Then $\rho_{\leq (n-1)}
\circ \jmath_2 \circ \jmath_1$ has image in $B^{n+1}(A,d) \subset
A^{n+1}$. Take a splitting of the map $d\colon A^{n} \surj
B^{n+1}(A,d)$, say
 $$
 \varrho\,\colon\, B^{n+1}(A,d) \,\too\, A^{n}\, ,
 $$
and
finally define
 $$
 \rho_n|_{N^n}\,=\,
 \varrho \circ \rho_{\leq (n-1)} \circ \jmath_2 \circ \jmath_1.
 $$

Now, define $d_n$ as follows. On $C^n$, we set
$d_n|_{C^n}=0$. On $N^n$, we put $d_n|_{N^n}= \jmath_3\circ
\jmath_2\circ \jmath_1$, where
 $$
 \jmath_3\,\colon\,  Z^{n+1}(\bV^{\leq (n-1)},d) \,\inc\,
 (\bV^{\leq (n-1)})^{n+1}\,= \,(\bbw\nolimits^{\geq 2} V^{\leq
 (n-1)})^{n+1}
 $$
is the inclusion. Clearly condition (ii) holds.

To check condition (iii), we need to verify that $\rho_{\leq (n-1)}
\circ d_n= d \circ \rho_n$. On $C^n$, we have $\rho_{\leq (n-1)}
\circ d_n=0$ and $d \circ \rho_n=d \circ \imath_3\circ
\imath_2\circ \imath_1=0$. On $N^n$, we have
 $$
 \rho_{\leq (n-1)} \circ d_n\,=\,
 \rho_{\leq (n-1)}\circ \jmath_3\circ \jmath_2\circ \jmath_1\,=\,
 d \circ \varrho \circ \rho_{\leq (n-1)} \circ \jmath_2 \circ
 \jmath_1\,=\,d\circ \rho_n
 $$
as $d \circ \varrho=\Id$. Therefore, condition (iii) holds.

Consider the inclusion $j\colon (\bV^{\leq (n-1)},d) \inc
(\bV^{\leq n},d)$ and the cokernel
 $$
 B=(\bV^{\leq n})/(\bV^{\leq (n-1)}) .
 $$
Then $(B,d)$ is a graded differential algebra, and $B^i=0$ for
all $i<n$, and also, $B^n=V^n= C^n\oplus N^n$. We have
$(\bV^{\leq
n})^{n+1} =(\bV^{\leq (n-1)})^{n+1}$ as $V^1=0$, hence
$B^{n+1}\,=\,0$. Therefore,
 $$
 j^*\,\colon\, H^k(\bV^{\leq (n-1)},d) \,
 \longrightarrow\, H^k(\bV^{\leq n},d)
 $$
is an isomorphism for all $k<n$. As $\rho_{\leq n}^*  \circ j^* =
\rho_{\leq (n-1)}^*$, we have that
 $$
 \rho_{\leq n}^*\,\colon\,
 H^k(\bV^{\leq n},d) \,\too\, H^k (A,d)
 $$
is an isomorphism for all $k<n$.

To deal with the cases where $k=n, n+1$, consider the long exact
sequence associated to $\bV^{\leq (n-1)} \inc \bV^{\leq n}\too B$,
 \begin{equation}\label{eqn:long}
 \begin{array}{rlclcl}
 0 \, \, \to & H^n(\bV^{\leq (n-1)},d) &\stackrel{j^*}{\too} &
 H^n(\bV^{\leq n},d)
 &{\to}&
 H^n(B,d)=B^n=C^n\oplus N^n\\[6pt]
 \stackrel{\bd^*}{\too} & H^{n+1}(\bV^{\leq (n-1)},d) &\to &
 H^{n+1}(\bV^{\leq n},d) &\to & 0 .
 \end{array}
 \end{equation}
For $x=u+w\in B^n=C^n\oplus N^n$, we have $\bd^*(x)=
[d(u+w)]=[\jmath_3 \circ\jmath_2 \circ\jmath_1(w)]=\jmath_1(w)$.

Therefore the exact sequence eqn. \eqref{eqn:long} splits into two
short exact sequences:
 \begin{equation}\label{eqn:short1}
 0\too H^n(\bV^{\leq (n-1)},d) {\too} H^n(\bV^{\leq n},d)
 \too C^n \too 0
 \end{equation}
and
 \begin{equation}\label{eqn:short2}
 0\too N^n \stackrel{\jmath_1}{\too} H^{n+1}(\bV^{\leq (n-1)},d) \too
 H^{n+1}(\bV^{\leq n},d) \too 0
 \end{equation}

{}From eqn. \eqref{eqn:short1}, we have
 $$
 \begin{array}{ccccccccc}
 0 &\too& H^n(\bV^{\leq (n-1)},d) &\too & H^n(\bV^{\leq n},d)
 & \too & C^n & \too&  0 \\
  && \| && {}_{\rho_{\leq n}^*}\downarrow && \downarrow\\
 0 &\too& H^n(\bV^{\leq (n-1)},d) &\stackrel{\rho_{\leq (n-1)}^*}{\too}
 & H^n(A,d)
 & \too & C^n & \too&  0
 \end{array}
 $$
We note that the right vertical arrow is the identity map. Indeed, it
sends $u\in
C^n$ to the class of $\rho_n(u)=\imath_3(\imath_2(\imath_1(u)))$
in $\coker(\rho_{\leq (n-1)}^*\colon H^n(\bV^{\leq (n-1)},d) \too
H^n(A,d))$, which is $u$ itself. Thus the middle vertical arrow is
an isomorphism, proving condition (iv) for $k=n$.

In eqn. \eqref{eqn:short2}, the
homomorphism $\jmath_1$ is the inclusion of
 $$
 N^n\,=\, \ker \left( \rho_{\leq (n-1)}^*\colon H^{n+1}(\bV^{\leq
 (n-1)},d) \too H^{n+1}(A,d)\right)
 $$
in $H^{n+1}(\bV^{\leq (n-1)},d)$. So $\rho_{\leq
(n-1)}^*$ induces an inclusion
 $$
 H^{n+1}(\bV^{\leq (n-1)},d) /N^n \,\cong\,
 H^{n+1}(\bV^{\leq n},d) \,\inc\, H^{n+1}(A,d)\, .
 $$
This map actually coincides with $\rho_{\leq n}^*$, since
$(\bV^{\leq (n-1)})^{n+1}=(\bV^{\leq n})^{n+1}$. This proves
that condition (v) holds.

\begin{remark} \label{rem:Nn}
{\rm Note that $d\colon N^n\too (\bV)^{n+1}$ is always injective,
and $d|_{C^n}=0$ for all $n$.}
\end{remark}

\begin{remark} \label{rem:extra}
{\rm If $(A,0)$ is a differential algebra with zero differential,
then the minimal model $\rho\colon (\bV,d) \too (A,0)$ constructed
above has the property that $\rho(N^n)=0$ for all $n$.}
\end{remark}

\section{Minimal models and $G$--actions} \label{sec:G-equivariant}

Let $G$ be a reductive complex Lie group. An \textit{action} of
the group $G$ on a differential algebra $(A,d)$ is a
representation $r\, :\, G\,\too\, \text{GL}(A)$ such that
\begin{itemize}
 \item $r(z)(A^n)\, =\, A^n$ for all $n\geq 0$ and $z\, \in\, G$,
 \item $r(z)(v_1\cdot v_2)\, =\, r(z)(v_1)\cdot r(z)(v_2)$ for all
  $v_1, v_2\, \in\, A$, and
 \item $r(z)(d v)\, =\, d(r(z)(v))$ for all $z\, \in\, G$ and $v\,
  \in\, A$.
\end{itemize}
If $G$ acts on $(A,d)$, then we say that $(A,d)$ is a
\emph{$G$--differential algebra}.

A \emph{$G$--minimal differential algebra} is a minimal
differential algebra $(\bV,d)$ on which $G$ acts satisfying the
condition that each graded vector space $V^n$,
$n\geq 0$, is preserved by the action of $G$. A \emph{$G$--minimal
model} of a $G$--differential algebra $(A,d)$ is a $G$--minimal
differential algebra $(\bV,d)$ such that there is a
$G$--equivariant map
 $$
 \rho\,\colon\, (\bV,d)\,\too\, (A,d)
 $$
which is a quasi--isomorphism.

Note that a $G$--minimal model is in
particular a minimal model.

\begin{prop} \label{prop:Gmm}
 Let $(A,d)$ be a $1$--connected $G$--differential algebra.
 Then there exists a $G$--minimal model $(\bV,d)$ of $(A,d)$.
\end{prop}

\begin{proof}
The construction of a minimal model
in Section \ref{sec:formal} works in the context
of $G$--differential algebras. All we need is to substitute the
vector spaces $V^n$ in the construction
by $G$--representations. We note that the reductivity of the group
$G$ ensures that any short exact sequence of $G$--modules splits.
\end{proof}

Let $(X\, ,x_0)$ be a one--pointed compact connected Riemann
surface of genus $g\geq 2$, and, as in Section
\ref{sec:introduction2}, let $\cN_X$ be the
moduli space of stable vector bundles over $X$ of
rank two and fixed determinant ${\mathcal O}_X(x_0)$.
Then the mapping class group acts on the
cohomology ring $H^*(\cN_X,\QQ)$, with the action factoring
through an action of $\Sp(2g,\ZZ)$; moreover, this action extends
to an action of $\Sp(2g,\CC)$ on $H^*(\cN_X,\CC)$ (Proposition
\ref{prop:action}). We will now show that a similar result holds
for the rational homotopy groups.

\begin{thm} \label{thm:main2}
The mapping class group acts on the homotopy groups
$\pi_*(\cN_X)$. The induced action on the rational homotopy groups
$\pi_*(\cN_X)\ox_\ZZ \QQ$ factors through an action of the
symplectic group ${\rm Sp}(2g,\ZZ)$. This action extends uniquely
to an action of ${\rm Sp}(2g,\CC)$ on $\pi_*(\cN_X)\ox_\ZZ \CC$.

Let $(\bV,d)$ be the ${\rm Sp}(2g,\CC)$--minimal model,
provided by Proposition \ref{prop:Gmm}, for the
$1$--connected ${\rm Sp}(2g,\CC)$--differential algebra
$(H^*(\cN_X,\CC)\, ,0)$. Then
  $$
  V^n\,\cong\, (\pi_n(\cN_X)\ox_\ZZ \CC)^*
  $$
as  ${\rm Sp}(2g,\CC)$--modules.
\end{thm}

\begin{proof}
First note that the formality of $\cN_X$ (Theorem \ref{prop:DGMS})
means that $(\bV\, ,d)$ is also the minimal model of $\cN_X$.

Now, let $\eta\in \G_g^1$, where $\G_g^1$ is
the mapping class group of
$(X\, , x_0)$. As we noted prior to Proposition \ref{prop:action},
the element $\eta$
acts on $\cN_X$ by a diffeomorphism $f_\eta \colon \cN_X\too \cN_X$.
Hence we have an action on the free homotopy groups of $\cN_X$. As
$\cN_X$ is simply connected, the free homotopy groups of $\cN_X$
coincide with the homotopy groups of $\cN_X$. So we have an
induced map
 $$
 \rho(\eta)\,\colon\, \pi_*(\cN_X)\,\too\, \pi_*(\cN_X)\, .
 $$

The diffeomorphism $f_\eta$ induces a map on differential forms,
 \begin{equation}\label{eqn:omega}
 f_\eta^*\,\colon\, (\Omega^*(\cN_X,\CC),d)\,\too\,
 (\Omega^*(\cN_X,\CC),d)\, ,
 \end{equation}
which lifts to a map on the minimal model
 \begin{equation}\label{eqn:w-h}
 \widehat{f}_\eta^*\,\colon\, (\bV,d)\,\too\, (\bV,d)\, .
 \end{equation}
Such a lift is not unique; it is only unique up to homotopy
of maps of differential algebras \cite{DGMS}. However, the induced map
on the indecomposables,
 \begin{equation}\label{eqn:o-b-f}
 \widetilde{f}_\eta^*\colon
  V=\frac{\bV}{\bigwedge^{\geq 2}V}\too V=\frac{\bV}{\bigwedge^{\geq
  2}V}\, ,
 \end{equation}
is unique \cite[Proposition 2.12]{VN}, and moreover, it coincides
with the dual of the map
 $$
 \rho(\eta)\ox\CC\,:\,\pi_*(\cN_X)\ox_\ZZ \CC\,\too\,
 \pi_*(\cN_X)\ox_\ZZ \CC\, ,
 $$
under the isomorphism of vector spaces $V^n\cong
(\pi_n(\cN_X)\ox_\ZZ\CC)^*$ (see \cite[page 259]{DGMS}).

Let
\begin{equation}\label{eqd.pi}
 \pi\,:\,\G_g^1\,\too\, \text{Aut}(H_1(X, {\mathbb Z}))
\, =\, \Sp(2g,\ZZ)
\end{equation}
be the natural projection of the
mapping class group onto the symplectic group. Then the automorphism
of cohomology
 $$
 \overline{f}_{\eta}^*\,\colon\, H^*(\cN_X,\CC)\,\too\,
 H^*(\cN_X,\CC)
 $$
induced by $f_\eta^*$ in eqn. \eqref{eqn:omega} coincides with the
action of $\pi(\eta)$ on the cohomology, where $\pi$ is the above
projection. The map $\widehat{f}_\eta^*$ in eqn. \eqref{eqn:w-h}
evidently induces the above automorphism $\overline{f}_{\eta}^*$.

The minimal model $(\bV,d)$ has an action of $\Sp(2g,\ZZ)$. Indeed,
by Proposition \ref{prop:Gmm}, the group $\Sp(2g,{\mathbb C})$ acts
on $(\bV,d)$ and this restricts to an action of
$\Sp(2g,\ZZ)\subset \Sp(2g,\CC)$. The homomorphism
$\overline{f}_{\eta}^*$ is induced by the action of
$\pi(\eta)$ on $(\bV,d)$, where $\pi$
is the projection in eqn. \eqref{eqd.pi}. Therefore the map
 $$
 \widetilde{f}_\eta^*\,\colon\, V^n \,\too\, V^n
 $$
defined in eqn. \eqref{eqn:o-b-f} coincides with the action of
$\pi(\eta)$ on $V^n$. Hence, under the isomorphism $V^n\, \cong
\, (\pi_n(\cN_X)\ox_\ZZ \CC)^*$, the actions of $\pi(\eta)$ and
$(\rho(\eta)\ox \CC)^*$ coincide.

If $\eta\, \in\, \Gamma^1_g$
belongs to the Torelli group, then $\rho(\eta)\ox \CC$ must be
the identity map of
$\pi_n(\cN_X)\ox_\ZZ \CC$, and hence $\rho(\eta)\ox
\QQ=\Id$ on $\pi_n(\cN_X)\ox_\ZZ\QQ$. This proves that the action
of the mapping class group on $\pi_*(\cN_X)\ox_\ZZ\QQ$ factors
through an action of $\Sp(2g,\ZZ)$. Moreover, this action
coincides with the restriction of the $\Sp(2g,\CC)$--action on
$V^n$ to the subgroup $\Sp(2g,\ZZ)\,\subset\, \Sp(2g,\CC)$ under
the isomorphism $(\pi_n(\cN_X)\ox_\ZZ \CC)^*\cong V^n$. So the
action of $\Sp(2g,\ZZ)$ on $\pi_*(\cN_X)\ox_\ZZ\CC$ extends to an
action of $\Sp(2g,\CC)$. Since $\Sp(2g,\ZZ)$ is Zariski dense in
$\Sp(2g,\CC)$ \cite{Bo}, the extension is unique. This completes
the proof of the theorem.
\end{proof}

Let $G\,=\,\Sp(2g,\CC)$, and let $(\bV,d)$ be the $G$--minimal
model of $(H^*(\cN_X,\CC),0)$. Then we may decompose $V^n$ into
irreducible $G$--representations. Let $\{\G_i\}_{i\in \Lambda}$ be
a complete set of irreducible $G$--representations, where
$\Lambda$ parametrizes the isomorphism classes of irreducible
$G$--representations. So
 $$
 V^n= \bigoplus_{i\in \Lambda} a_{i,n} \G_i
 $$
for some set of integers $a_{i,n}\geq 0$.

\section{The minimal model of $\cN_X$ for $g=2$} \label{sec:moduli}

In this section we will assume that $X$ is a compact connected
Riemann surface of genus two. In this case,
the moduli space $\cN_X$, whose dimension is now three, can
be described explicitly \cite{Ne1, NR}. It turns out to be
isomorphic to the intersection of two quadrics in $\PP^5$.

The integral cohomology ring of $\cN_X$ has no torsion
\cite[Section 10]{Ne0}. Let $h\in H^2(\cN_X,\ZZ)$ be the
hyperplane class. Note that, by the Lefschetz hyperplane theorem,
$H^2(\cN_X,\ZZ)\,=\, H^2(\PP^5,\ZZ)\,=\,\ZZ$. So $h\,=\,\alpha$,
the generator of the ample cone. The intersection of two quadrics
in $\PP^5$ contains many lines $\PP^1\subset \cN_X\subset \PP^5$.
Let $l\,\in\, H^4(\cN_X,\ZZ)$ be the Poincar{\'e} dual of such a line.
Then $h\cup l =[\cN_X]$, so $H^4(\cN_X,\ZZ)\cong \ZZ$ is generated
by $l$. We have $h^3=h\cup h \cup h = 4[\cN_X]$, as the degree of
$\cN_X \subset\PP^5$ is $4$. Therefore, we conclude that $h\cup h=
4l$. Finally, $H^3(\cN_X,\ZZ)\cong H^1(X,\ZZ)^*$, so
$H^3(\cN_X,\ZZ)\cong W_0$, the standard
$\Sp(4,\ZZ)$--representation $W_0=\ZZ^{4}$. Moreover, the pairing
 $$
 H^3(\cN_X,\ZZ)\ox_\ZZ H^3(\cN_X,\ZZ) \too H^6(\cN_X,\ZZ)\cong \ZZ
 $$
is perfect (Poincar{\'e} duality) and $\Sp(4,\ZZ)$--equivariant, so it
is equivalent to the standard symplectic form on $W_0$. The
conclusion is that
  $$
  \left\{
  \begin{array}{lll}
  H^0(\cN_X,\ZZ)&=&  \la 1\ra\, ,  \\
  H^1(\cN_X,\ZZ)&=&0\, , \\
  H^2(\cN_X,\ZZ)&=&\la h\ra\, , \\
  H^3(\cN_X,\ZZ)& \cong & W_0\, , \\
  H^4(\cN_X,\ZZ)&=&\la l \ra\, , \\
  H^5(\cN_X,\ZZ)&=&0\, ,\\
  H^6(\cN_X,\ZZ)&=&\la [\cN_X]\ra\, .
  \end{array}
  \right.
  $$

This can also be seen by using Proposition \ref{prop:cohomology},
at least for rational coefficients. Since $I_1=(\a,\b,\g)$ and
$I_2=(\a^2+\b, \a\b+\g,\a\g)$, Proposition \ref{prop:cohomology}
says that
 $$
 H^*(\cN_X,\QQ)\,=\, \frac{\QQ [\a, \b, \g]}{I_2} \oplus \left(W
 \ox \frac{\QQ [\a, \b, \g]}{I_1}\right) \,\cong\,
 \frac{\QQ[\a]}{(\a^4)}\oplus W\, ,
 $$
where $\b= -\a^2$ and $\g=\a^3$ in this ring, and $\g_i\cup\g_j= -
\frac{1}{4}(\g_i\cdot \g_j) \a^3$, for any $\g_i, \g_j\in W$. Note
that $\b=-4l$. Note that $W=W_0\ox_\ZZ \QQ$ is the standard
$\QQ$--representation of $\Sp(4,\ZZ)$.

Now we pass on to compute the minimal model $(\bV,d)$ of $\cN_X$
by computing the minimal model of its cohomology algebra
$H^*(\cN_X,\CC)$. This is possible because $\cN_X$ is formal by
Theorem \ref{prop:DGMS}. By Proposition \ref{prop:Gmm},
$(\bV,d)$ is a $G$--minimal model for $G\,=\,\Sp(4,\CC)$.

The irreducible representations of $\Sp(4,\CC)$ are labeled by
pairs $(a,b)$, $a,b\geq 0$, such that the corresponding
representation $\G_{a,b}$ has highest weight
$aL_1+b(L_1+L_2)=(a+b)L_1+bL_2$, where $L_1$ and $L_2$ are the
orthogonal generators (with respect to the Killing form) of the
weight lattice; see \cite[Part III, Section 16]{FH}.

The standard representation $W_c=W\ox_\QQ \CC=\CC^4$ of
$\Sp(4,\CC)$ is $W_c=\G_{1,0}$, whereas the irreducible
$\Sp(4,\CC)$--representation $\bw^2_0 W_c$ is $\G_{0,1}$. Some
easy cases are dealt with in \cite[Part III, Section 16]{FH},
  $$
  \left\{
  \begin{array}{l}
  \bw^2 \G_{1,0} = \bw^2_0 W_c \oplus \CC= \G_{0,1} \oplus\G_{0,0}\, ,\\
  \Sym^a \G_{1,0} = \G_{a,0}\, , \\
  \G_{0,1}\ox \G_{1,0} = W_c\ox \bw^2_0 W_c = \G_{1,1}\oplus
  \G_{1,0}\, .
  \end{array} \right.
  $$

We define a partial order in the set of weights of $\Sp(4,\CC)$ as
follows:
  $$
  (a,b) \leq (c,d) \Leftrightarrow \left\{
  \begin{array}{l} a+b \leq c+d,
  \\ a+2b \leq c+2d. \end{array} \right.
  $$
This corresponds to the fact that the weights of the
representation $\G_{a,b}$ are a subset of the convex hull of the
weights of $\G_{c,d}$. Otherwise said, $(a,b)\geq 0$ means that
the highest weight $(a+b)L_1+bL_2$ is a linear combination with
non-negative coefficients of the positive roots (see \cite{FH}).
(We point out that this order is defined in \cite[page 47]{Hum}
with the difference that in \cite{Hum}, $(a,b)\geq 0$ means that
$(a+b)L_1+bL_2$ is a linear combination with non-negative
\textit{integer} coefficients of the positive roots. This is
equivalent to $a+b\geq 0$, $a+2b\geq 0$ and $a+2b\equiv 0\pmod
2$.)

In particular, for representations $\G_{a_1,b_1}$ and
$\G_{a_2,b_2}$, the sub--representations $\G_{c,d}$ of the tensor
product $\G_{a_1,b_1} \ox \G_{a_2,b_2}$ satisfy the condition
  $$
  (c,d) \leq (a_1,b_1)+(a_2,b_2)=(a_1+a_2,b_1+b_2)\, ,
  $$
and furthermore, there is exactly one sub--representation
(the Cartan component) satisfying the
equality. Note that this says in particular that $\G_{a,b}\subset
W_c^{\ox a} \ox (\bw^2_0 W_c)^{\ox b}$ appears with multiplicity
one.

We compute the $\Sp(4,\CC)$--minimal model $(\bV,d)$ following the
mechanism laid out in Section \ref{sec:formal} and Proposition
\ref{prop:Gmm}.

\begin{prop} \label{prop:g=2}
Let $(\bV,d)$ be the minimal model of $\cN_X$ for a curve $X$ of genus
$g=2$. Then we have, as $\Sp(4,\ZZ)$--representations,
  $$
  \left\{ \begin{array}{l}
   V^2= \G_{0,0} \\
   V^3= \G_{1,0}\\
   V^4= \G_{1,0}\\
   V^5= \G_{0,1} \oplus \G_{0,0}\\
   V^6= \G_{2,0}\oplus \G_{0,1} \\
   V^7= \G_{1,1} \oplus \G_{2,0}\oplus \G_{1,0}\, .
     \end{array} \right. \hspace{8cm}
  $$
\end{prop}

\begin{proof}
 Abbreviating $H^*(\cN_X)$ for $H^*(\cN_X,\CC)$, we have
  $$
  V^2=C^2=H^2(\cN_X)=\la h \ra \cong \G_{0,0}\,.
  $$
 Recall that $d|_{C^n}=0$, for any $n$.

 In the next step, we have $V^3=C^3\oplus N^3$, with
  \begin{eqnarray*}
  C^3 &=& \coker \left(H^3(\bV^{\leq 2})=0 \too H^3(\cN_X)\right)= H^3(\cN_X)
  \cong W_c\cong \G_{1,0}\, ,\\
  N^3 &=&\ker \left(H^4(\bV^{\leq 2})=\la h^2\ra \too H^4(\cN_X)=\la
  h^2\ra\right)=0.
  \end{eqnarray*}

 For $n=4$, we have $V^4=C^4\oplus N^4$, with
  \begin{eqnarray*}
  C^4 &=& \coker \left(H^4(\bV^{\leq 3})=\la h^2\ra \too H^4(\cN_X)\la
  h^2\ra \right) = 0 ,\\
  N^4 &=&\ker \left(H^5(\bV^{\leq 3})
  \too H^5(\cN_X)=0\right)=H^5(\bV^{\leq 3}) \\
  &=& V^3\cdot
  V^2
  \cong \G_{1,0}\ox \G_{0,0}=\G_{1,0}.
  \end{eqnarray*}
 The differential $d\colon N^4\too V^3\cdot V^2 \subset \bV$ is an
 isomorphism.

 We continue with $V^5=C^5\oplus N^5$, where
  \begin{eqnarray*}
  C^5 &=& \coker \left(H^5(\bV^{\leq 4}) \too H^5(\cN_X)=0 \right)=0\, ,\\
  N^5 &=& \ker \left(H^6(\bV^{\leq 4})=\bw^2 V^3 \oplus \la h^3\ra
  \too H^6(\cN_X)=\la h^3\ra\right) \\ & \cong & \bw^2 V^3 \cong \bw^2
  \G_{1,0}= \G_{0,1}\oplus \G_{0,0} \, .
  \end{eqnarray*}
The differential $d\,\colon\, N^5\,\too\, \bw^2 V^3 \oplus \la
h^3\ra$ is an isomorphism of $N^5$ with the kernel of the map
$\bw^2 V^3 \oplus \la h^3\ra\too \la h^3\ra$. This map is the sum
of a multiple of the intersection product $\bw^2 V^3 \too \CC
\cong \la h^3\ra$ in the first summand, and the identity in the
second summand.

For $n=6$, we have $V^6=C^6\oplus N^6$. Now
 $$
 C^6 \,=\, \coker
 \left(H^6(\bV^{\leq 5}) \surj H^6(\cN_X)=\la h^3\ra \right)=0\, ,
 $$
since $h^3\in H^6(\bV^{\leq 5})$. Moreover $C^k=0$ for $k>6$ since
$H^k(\cN_X)=0$. Also for all $k\geq 6$, we have
$N^k=H^{k+1}(\bV^{\leq (k-1)})$, since $H^{k+1}(\cN_X)=0$. Now
  \begin{eqnarray*}
  (\bV^{\leq 5})^6 &=& (V^3 \cdot V^3) \oplus
 (V^2\cdot V^2\cdot V^2) \oplus (V^4 \cdot V^2)\, , \\
  (\bV^{\leq 5})^7 &=& (V^4 \cdot V^3) \oplus
 (V^3\cdot V^2\cdot V^2) \oplus (V^5 \cdot  V^2)\, .
  \end{eqnarray*}
The space of coboundaries is $B^7(\bV^{\leq 5})\,=\, d(V^4\cdot
V^2)\,=\,V^3\cdot V^2\cdot V^2$. The differential $d$ maps $(V^4
\cdot V^3) \oplus (V^5 \cdot V^2)$ onto $\bw^2 V^3 \cdot V^2
\oplus \la h^4\ra$, and it has kernel isomorphic to $\ker
(V^4\cdot V^3 \too \la h^4\ra)$. But
  $$
  V^4\cdot V^3=V^4\ox V^3\cong \G_{1,0}\ox\G_{1,0}=\Sym^2\G_{1,0}\oplus
  \bw^2\G_{1,0}\cong \G_{2,0}
  \oplus \G_{0,1}\oplus \G_{0,0}\, .
  $$
So the conclusion is that
  $$
  N^6= H^7(\bV^{\leq 5})=\frac{Z^7(\bV^{\leq 5})}{B^7(\bV^{\leq
  5})} \cong\G_{2,0} \oplus \G_{0,1}\, ,
  $$
and the differential $d\colon N^6 \too (V^4 \cdot V^3) \oplus (V^5
\cdot V^2)$ is the sum of the two maps $d\colon N^6=\G_{2,0}
\oplus \G_{0,1} \too V^4 \cdot V^3=\G_{2,0} \oplus \G_{0,1}\oplus
\G_{0,0}$ which is injective, and $d\colon N^6=\G_{2,0} \oplus
\G_{0,1} \too V^5 \cdot V^2=\G_{0,1}\oplus \G_{0,0}$ mapping onto
the $\G_{0,1}$ summand.

 The next case is $V^7=C^7\oplus N^7=N^7$. Then
  $$
  \begin{array}{l} \\[-8pt]
(\bV^{\leq 6})^7 = (V^5 \cdot V^2)\oplus
(V^4 \cdot V^3) \oplus (V^3\cdot V^2\cdot V^2), \\[6pt]
(\bV^{\leq 6})^8 = (V^6\cdot V^2) \oplus (V^5 \cdot V^3) \oplus
(V^4\cdot V^4) \oplus (V^4\cdot V^2\cdot V^2) \oplus (V^3 \cdot
V^3\cdot V^2)\oplus \la h^4\ra .\\
  \end{array}
  $$
The space of coboundaries is $B^8(\bV^{\leq 6}) =(V^3 \cdot
V^3\cdot V^2) \oplus \la h^4\ra$, from our knowledge of $d$ on
both $V^5$ and $V^4$. To compute
  \begin{eqnarray*}
   N^7 &=& H^8(\bV^{\leq 7})=\frac{Z^8(\bV^{\leq 6})}{B^8(\bV^{\leq
   6})} = \\
  &=& \ker \left(d\colon (V^6\cdot V^2)
  \oplus (V^5 \cdot V^3) \oplus (V^4\cdot V^4) \oplus
   (V^4\cdot V^2\cdot V^2) \too \bV\right)\, ,
  \end{eqnarray*}
  we look at each summand,
  \begin{eqnarray*}
  && d\colon V^5 \cdot V^3\too \bw^3 V^3 \oplus \, V^3\cdot (V^2)^3, \\
  && d\colon  V^4\cdot V^2\cdot V^2\isom V^3 \cdot (V^2)^3, \\
  && d\colon  V^6\cdot V^2 \inc (V^4 \cdot V^3\cdot V^2) \oplus (V^5 \cdot
  V^2\cdot V^2) \\
  && d \colon  V^4\cdot V^4 \inc V^4 \cdot V^3\cdot V^2.
  \end{eqnarray*}
  So
  $N^7=K_1\oplus K_2$, where $K_1=\ker((V^5 \cdot V^3) \oplus
   (V^4\cdot V^2\cdot V^2) \too \bV)$ and
   $K_2=\ker( (V^6\cdot V^2) \oplus (V^4\cdot
  V^4)\too \bV)$. Clearly, $K_1\cong \ker
   (V^5\ox V^3 \too \bw^3 V^3)$, but
  $V^5\ox V^3\cong (\G_{0,1} \oplus \G_{0,0})\ox
  \G_{1,0} =
  \G_{1,1}\oplus \G_{1,0} \oplus\G_{1,0}$
   and $\bw^3 V^3\cong V^3 \cong \G_{1,0}$, so $K_1\cong
   \G_{1,1}\oplus \G_{1,0}$. On the other hand, $d$
   maps $V^4\cdot V^4=\Sym^2 V^4 \cong \G_{2,0}$ to the
   corresponding summand in $V^4 \cdot V^3\cdot V^2\cong
   \G_{1,0}\ox \G_{1,0}=\G_{2,0}\oplus \G_{0,1}\oplus \G_{0,0}$,
   and $d$ maps $V^6\cdot V^2 \cong \G_{2,0}\oplus \G_{0,1}$
   injectively to $(V^4 \cdot V^3\cdot V^2) \oplus (V^5 \cdot
  V^2\cdot V^2) \cong (\G_{2,0}\oplus \G_{0,1}\oplus \G_{0,0})
  \oplus (\G_{0,1}\oplus \G_{0,0})$. Thus $K_2\cong \G_{2,0}$.
  This concludes that $N^7=\G_{1,1}\oplus \G_{2,0} \oplus
  \G_{1,0}$, and the proof of the proposition is complete.
\end{proof}

We may carry on the process as long as we want, but the
calculations get more involved, since we must keep track of the
irreducible summands of $\bV$ onto which $d|_{V_n}$ maps for each
$n$. It is easier to find the ``leading representation''. We need
a preliminary result.

\begin{lemma} \label{lem:n(a,b)}
 For any ${\rm Sp}(4,\CC)$--irreducible representation $\G_{a,b}$,
$a,b\geq 0$,
 if $\G_{a,b}\subset V^n$, then $n\geq n(a,b)$, where
   $$
   n(a,b)=\left\{ \begin{array}{ll} 2a+4b+1, \qquad & \text{if $b\geq
   1$ or $(a,b)=(1,0)$,} \\
   2a+2, & \text{if $b=0$ and $a\neq 1$.}
  \end{array}
   \right.
   $$
\end{lemma}

\begin{proof}
 We shall prove this by induction on $n$. By Proposition \ref{prop:g=2},
the result is true for $1\leq n \leq 7$. So suppose $n\geq 8$.
 Let $U\subset V^n$ be a sub--representation with $U\cong\G_{a,b}$.
 We want to
 prove that $n(a,b)\leq n$, so we may assume that $(a,b)\neq (0,0)$,
 $(1,0)$, $(0,1)$, $(1,1)$ and $(2,0)$. As
$H^n(\cN_X)=0$, we have that $C^n=0$. So
 $V^n=N^n$, and in particular $d\colon V^n
 \too (\bV)^{n+1}$ is injective.
 Hence
  $$
  U \cong d(U)\subset d(V^n) \subset (\bbw\nolimits^{\geq 2}
  V)^{n+1} = \sum_{n_1+\ldots +n_r=n+1, \ r\geq 2} V^{n_1}\cdots
  V^{n_r}
  \, .
  $$
 The projection of $d(U)$ to some of these summands
must be non--zero. Hence there
 exists
  $$
  U' \subset V^{n_1}\cdots
  V^{n_r} \subset V^{n_1}\ox \cdots \ox V^{n_r},
  $$
 for some $r\geq 2$, with $n_1+\ldots + n_r=n+1$, $U'\cong U\cong \G_{a,b}$.
 Note that all $n_i<n$, because $n_i\geq
 2$, $1\leq i \leq r$.
 Decomposing each $V^{n_i}$ into $\Sp(4,\CC)$--irreducible
 representations, there must be $(a_1,b_1),\ldots, (a_r,b_r)$
 such that
 \begin{equation}\label{eqn:inclusion}
  \G_{a,b}\subset \G_{a_1,b_1}\ox \cdots \ox
  \G_{a_r,b_r}\, ,
 \end{equation}
 with $\G_{a_i,b_i} \subset V^{n_i}$. Applying the
 induction hypothesis it follows that
 $n_i\geq n(a_i,b_i)$ for all $1\leq i\leq r$.
 We note that eqn. \eqref{eqn:inclusion} implies that
 $(a,b)\leq (a_1,b_1) +\ldots +(a_r,b_r)$, that is
  \begin{equation} \label{eqn:inequ}
  \left\{ \begin{array}{l}
  a+b \leq \sum (a_i+b_i), \\
  a+2b \leq \sum (a_i+2b_i).
  \end{array}\right.
  \end{equation}
 If $a,b\geq 1$, we have
  \begin{equation} \label{eqn:claim}
  \begin{array}{lcl}
  n+1 &=& n_1+\ldots + n_r  \\
  &\geq& n(a_1,b_1) +\ldots +n(a_r,b_r) \\
  &\geq & \sum_{i=1}^r (2a_i+4b_i+1) \\
  &\geq& 2a+4b +r \\ &\geq & 2a+4b+2 \\ &=& n(a,b)+1 ,
  \end{array}
  \end{equation}
using eqn. \eqref{eqn:inequ}. So $n\geq n(a,b)$ in this case.

If $b=0$, $a\geq 3$, then eqn. \eqref{eqn:claim} proves that $n+1\geq
n(a,b)$. If there is equality, then $r=2$,
$n_i=n(a_i,b_i)=2a_i+4b_i+1$ and $a=a+2b=\sum (a_i+2b_i)$, for all
$i$. Since $a=a+b\leq \sum (a_i+b_i)$, we get $\sum(a_i+2b_i) =a
\leq \sum (a_i+b_i)$, so $b_i=0$ for all $i$. As also
$n_i=n(a_i,0)=2a_i+1$, we have that $a_i=1$, $n_i=3$. But then
$a=2$ which is a case treated before.
\end{proof}

\begin{thm} \label{thm:leading}
 Let $n\,\geq\, 4$. The decomposition of $V^n$ into
a direct sum of irreducible $\Sp(4,\CC)$--representations
 is as follows:
 \begin{itemize}
 \item[(i)] If $n=2m$ is even, then
  $$
  V^n= \G_{m-1,0} \oplus  \left(\bigoplus_{(a,b)<(m-1,0)}
  n_{ab}\G_{a,b}\right) ,
  $$
 with $n_{ab}\geq 0$.
 \item[(ii)] If $n=2m+1$ is odd, then
  $$
  V^n=
  \G_{m-2,1} \oplus \left(\bigoplus_{(a,b)<(m-2,1)}
  n_{ab}\G_{a,b}\right),
  $$
 with $n_{ab}\geq 0$.
 \end{itemize}
\end{thm}

\begin{proof}
(i) Let $n=2m\geq 4$. By Lemma \ref{lem:n(a,b)}, if
$\G_{a,b}\subset
 V^n$ then $2a+4b+1\leq n=2m$. This implies that
 $$
 (a,b)\leq (m-1,0),
 $$
 since $a+2b\leq m-1$ and $a+b\leq a+2b\leq m-1$.
So the leading representation in $V^n$ is $\G_{m-1,0}$.
We will show that it actually appears and that its multiplicity
is one.

 To see that there is $\G_{m-1,0}
 \subset V^n$, we shall prove by induction on $m$ that there is a
 sub--representation $U_{m-1}\subset V^n$, $U_{m-1}\cong
 \G_{m-1,0}$, such that $d(U_{m-1})\subset V^{n-2}\cdot V^3$.
 By Proposition \ref{prop:g=2}, this
 is true for $m=2,3$. Assume that it is true for $m-1\geq 3$ and
 let us prove it for $m$. So $U_{m-2}\subset V^{2m-2}$ and
 $d(U_{m-2})\subset V^{2m-4}\cdot V^3$. Then
 $d\colon  U_{m-2}\cdot V^3  \too V^{2m-4}\cdot V^3 \cdot V^3 \subset
 \bV$. But
 $$
 U_{m-2}\cdot V^3\cong  \G_{m-2,0}\ox \G_{1,0}
 $$
 contains a sub--representation $\widetilde{U}_{m-1}\subset  U_{m-2}\cdot
 V^3$ such that $\widetilde{U}_{m-1} \cong \G_{m-1,0}$. On the other hand,
  $$
  V^{2m-4}\cdot V^3 \cdot V^3= V^{2m-4}\ox \bw^2 V^3 \cong
  V^{2m-4}\ox \G_{0,1} \, .
  $$
 Decomposing $V^{2m-4}$ into irreducible representations
 $\G_{c,d}$, and noting that $(c,d)\leq (m-3,0)$ by induction hypothesis,
we see that if $\G_{a,b} \subset V^{2m-4}\cdot V^3 \cdot V^3$ then
 $\G_{a,b} \subset \G_{c,d}\ox \G_{0,1}$ for some $(c,d)\leq (m-3,0)$.
 Thus
  $$
  (a,b) \leq (c,d)+ (0,1)=(c,d+1) \leq (m-3,1)< (m-1,0).
  $$
 As a consequence,  $\G_{m-1,0}\not\subset V^{2m-4}\cdot V^3 \cdot
 V^3$, and so
 $d(\widetilde{U}_{m-1})=0$. This implies that $\widetilde{U}_{m-1}
 \subset Z^{n+1}(\bV,d) =B^{n+1}(\bV,d)$, since
 $H^{n+1}(\bV,d)=0$. There must exist a sub--representation
 $U_{m-1}\subset
 (\bV)^n$ with $d(U_{m-1}) = \widetilde{U}_{m-1}$. As $d$ maps
 $(\bbw^{\geq i} V)^{n} \too (\bbw^{\geq (i+1)} V)^{n+1}$, it cannot
 be $U_{m-1}\subset(\bbw^2 V)^{n}$, so the projection of $U_{m-1}$
 by $p\colon (\bw V)^n \too V^n$ is a sub--representation isomorphic to
 $\G_{m-1,0}$. (Here we are allowed to substitute $V^n$ by $U_{m-1} \oplus
 p(U_{m-1})^{\perp}$, where $p(U_{m-1})^{\perp}$ is a
 $\Sp(4,\CC)$--invariant
complement of $p(U_{m-1})\,\subset\, V^n$; this yields an
isomorphic minimal model and ensures that $d(U_{m-1}) \subset
V^{2n-2}\cdot V^3$).

Now let us compute the multiplicity of $\G_{m-1,0}$ in $V^n$. The
argument of the proof of Lemma \ref{lem:n(a,b)} implies that the
multiplicity of $\G_{m-1,0}$ in $V^n$ is at most the sum of the
multiplicities of $\G_{m-1,0}$ in $V^{n_1}\cdots V^{n_r}$, for the
different possibilities $n_1+ \ldots+n_r= n+1$, $r\geq 2$. Let
$(a,b)=(m-1,0)$. As in the proof of Lemma \ref{lem:n(a,b)}, for
any sub--representation $\G_{m-1,0}$ there are $(a_1,b_1), \ldots,
(a_r,b_r)$ such that
  \begin{equation} \label{eqn:claim2}
  \begin{array}{lcl}
  2m+1 &=& n+1 = n_1+\ldots + n_r   \\
  &\geq& n(a_1,b_1) +\ldots +n(a_r,b_r) \\
  &\geq& \sum_{i=1}^r (2a_i+4b_i+1) \\
  &\geq& 2a+4b +r = 2m-2 +r .
  \end{array}
  \end{equation}
In particular $r\leq 3$. If $r=3$, then
$n_i=n(a_i,b_i)=2a_i+4b_i+1$ for all $i$, and $a=a+2b=\sum
(a_i+2b_i)$. Since $a=a+b\leq \sum (a_i+b_i)$, we get $b_i=0$ for
all $i$. This implies that $a_i=1$ and $n_i=3$. But then $a=m-1=3$
and
$$
\G_{3,0}\not\subset V^3\cdot V^3\cdot V^3=\bw^3V^3 \cong \bw^3
W_c \cong W_c \cong \G_{1,0}\, .
$$

If $r=2$, then $2a+4b +1 \geq \sum (2a_i+4b_i) \geq 2a+4b$. So
$\sum(a_i+2b_i)=a+2b=a$. As before, this implies that $b_i=0$ for
all $i$. At most one of the $a_i$'s is bigger than $1$, so we can
put $(a_1,b_1)=(m-2,0)$, $(a_2,b_2)=(1,0)$. This corresponds to
the summand $\G_{m-2,0}\ox \G_{1,0} \subset V^{2m-4}\cdot V^3$.
This representation contains $\G_{m-1,0}$ with multiplicity one.

Since we know that the multiplicity of $\G_{m-1,0}$ in $V^n$ is
non--zero, we conclude that it is exactly one.

(ii) Let $n=2m+1\geq 5$. By Proposition \ref{prop:g=2}, the result
holds for $m=2,3$, so assume that $m\geq 4$.

 If $\G_{a,b}\subset V^n$, then by
 Lemma \ref{lem:n(a,b)}, we have that $2a+4b+1 \leq n=2m+1$, so
 $a+2b\leq m$. This implies that
 $$
 (a,b)\leq (m-2,1),
 $$
 since if $b\geq 1$ then $a+b \leq a+2b-1 \leq m-1$; and if
 $b=0$ then Lemma \ref{lem:n(a,b)} says that $2a+2 \leq n=2m+1$, so $a+b
 = a\leq m-1$. So the leading representation in $V^n$ is
 $\G_{m-2,1}$. We will show that it actually appears with
 multiplicity one.

As in the previous case, one can see using induction on $m$ that
there is a sub--representation $U_{m-1} \subset V^{2m+1}$,
with $U_{m-1}\cong
\G_{m-2,1}$, such that $d(U_{m-1})\subset V^{2m-1}\cdot V^3$.

To compute the multiplicity of $\G_{m-2,1}$ in $V^n$, let us find
the multiplicity of $\G_{m-2,1}$ in $V^{n_1}\cdots V^{n_r}$, for
$n_1+ \ldots+n_r= n+1$, $r\geq 2$. As $n=2m+1=n(m-2,1)$, there
must be equality in eqn. \eqref{eqn:claim} for $(a,b)=(m-2,1)$, which
means that $r=2$, $n_i=n(a_i,b_i)=2a_i+4b_i+1$ and $\sum
(a_i+2b_i)=a+2b=m$. Since $m-1=a+b\leq \sum (a_i+b_i)$, we have
$\sum (a_i+2b_i)  =m \leq \sum (a_i+b_i)$ +1 , so $\sum b_i\leq
1$. As least one $b_i$ is zero, say $b_2=0$. Then $a_2=1$,
$n_2=3$. Also $m-1 \leq \sum (a_i+b_i) \leq a_1 +2$ and
$m=a_1+2b_1+1$, implying that $(a_1,b)=(m-3,1)$ or $(m-1,0)$ and
$n_1=2a_1+4b_1+1=2m-1$. By induction hypothesis,
$\G_{m-1,0}\not\subset V^{2m-1}$, so the second case is ruled out.
The multiplicity of $\G_{m-2,1}$ in $\G_{m-3,1}\ox \G_{1,0}\subset
V^{2m-1}\cdot V^3$ is $1$. This proves that the multiplicity of
$\G_{m-2,1}$ in $V^n$ is one.
\end{proof}

\section{Sub--representations in the minimal model of $\cN_X$ for $g>2$}
\label{sec:general}

Suppose now that $X$ is a smooth irreducible projective complex
curve of genus $g\,>\, 2$. The action of $\Sp(2g,\CC)$ on the
cohomology algebra $H^*(\cN_X,\CC)$ of the moduli space $\cN_X$
gives an action of $\Sp(2g,\CC)$ on the minimal model $(\bV,d)$ of
$\cN_X$, by Proposition \ref{prop:Gmm}. By Theorem
\ref{thm:main2}, the action of $\Sp(2g,\CC)$ on the minimal model
$(\bV,d)$ is compatible with the action of $\Sp(2g,\CC)$ on the
complex homotopy groups $\pi_*(\cN_X)\ox_\ZZ \CC$.

The isomorphism classes of irreducible
$\Sp(2g,\CC)$--representations are labeled by $g$--tuples
$(a_1,\ldots, a_g)\in (\ZZ_{\geq 0})^g$ (see \cite[Part III,
Section 17]{FH}). The representation corresponding to
$(a_1,\ldots, a_g)$ is denoted by
  $$
\G_{(a_1,\ldots,a_g)} = \G_{a_1{\mathbf{e}}_1+\ldots +a_g{\mathbf{e}}_g},
  $$
where $\mathbf{e}_i=(0,\ldots,1,\ldots,0)$, with $1$ in the
$i$--th position and $0$ elsewhere. The $\Sp(2g,\CC)$--module
$\G_{(a_1,\ldots,a_g)}$ is characterized by its highest weight
$(a_1+a_2+\ldots +a_g) L_1+ (a_2+\ldots +a_g) L_2+ \ldots
+a_gL_g$, where $\{L_1,\ldots,L_g\}$ is the standard basis for the
weight lattice.

Let $W_c=\CC^{2g}$ be the standard representation of
$\Sp(2g,\CC)$. Then $W_c=\G_{{\mathbf{e}}_1}$, and
$\G_{{\mathbf{e}}_k}=\bw^k_0 W_c$ is the complexification of the
representation $\bw^k_0 W$ introduced in Section
\ref{sec:cohomology}.

We shall use two well--known facts: (1) the representation
$\G_{(a_1,\ldots,a_g)} \ox \G_{(b_1,\ldots,b_g)}$ contains
$\G_{(a_1+b_1,\ldots,a_g+b_g)}$ (actually this is the highest
weight representation appearing with multiplicity one); and (2)
the representation $\G_{(k-2,0,\ldots,0)}\ox \G_{(0,1,\ldots,0)}$
does not contain $\G_{(k,0,\ldots,0)}$ (this holds because the
weight $kL_1$ does not appear in the tensor product), and the
representation $\G_{(k-2,1,0,\ldots,0)}\ox \G_{(0,1,0,\ldots,0)}$
does not contain $\G_{(k,1,0,\ldots,0)}$.

\begin{thm} \label{thm:g>2}
Let $X$ be a complex smooth projective irreducible curve of genus
$g> 2$. Let $(\bV,d)$ be the minimal model of the moduli space
$\cN_X$. Then, as ${\rm Sp}(2g,\CC)$--representations, we have
   $$
   \left\{ \begin{array}{ll}
   V^2= \G_{0} \, ,\\
   V^3= \G_{{\mathbf{e}}_1}\, ,\\
   V^4= \G_{0}\, ,\\
   V^n= 0, \qquad & 5\leq n\leq 2g-2, \\
   V^{2g-1}= \G_{0} \, , \\
   V^{2g}= \G_{{\mathbf{e}}_1}\, ,\\
   V^{2g+1}= \G_{{\mathbf{e}}_2}\oplus \G_0\, .
   \end{array} \right. \hspace{5cm}
   $$
Moreover, for $n\geq 2g+2$, we have the following:
  \begin{itemize}
  \item[(i)] If $n=2(g+k-1)$ with $k\geq 2$, then $V^n$ contains
   $\G_{k {\mathbf{e}}_1}$.
  \item[(ii)] If $n=2(g+k)+1$ with $k\geq 1$, then $V^n$ contains
   $\G_{k {\mathbf{e}}_1+{\mathbf{e}}_2}$.
  \end{itemize}
\end{thm}

\begin{proof}
Clearly, $V^2\,=\,C^2\,=\,\la \a \ra\,\cong\, \G_0$,
 $$
 V^3\,=\,C^3\,=\,H^3(\cN_X)\,=\,W_c\cong \G_{{\mathbf{e}}_1}
 $$
and $V^4\,=\,C^4\,=\,\la \b\ra\,\cong\, \G_0$. Now
 $$
 \bV^{\leq 4}\,=\, \bigwedge (\a, \g_1,\ldots,\g_{2g},\b)\,=\,\AA_c\, ,
 $$
where $\AA_c=\AA\ox_\QQ\CC$ is the complexification of the
rational vector space defined in Lemma \ref{lem:kerF}. So the natural
homomorphism $\bV^{\leq 4}\,\longrightarrow\, H^*(\cN_X)$
is surjective. This implies that $C^n=0$ and
   $$
   V^n= N^n= \ker\left(H^{n+1}(\bV^{<n}) \too H^{n+1}(\cN_X)\right)
   $$
for all $n>4$.
Since $F_c:\AA_c\too H^*(\cN_X)$ is the complexification of the
map $F$ in eqn. \eqref{eqn:epim}, its kernel, ${\rm kernel} (F_c)$, has
the lowest degree element $q_g^1$, which is of degree $2g$. So
$V^n=N^n=0$ for all $5\leq n\leq 2g-2$. For $n=2g-1$, we have
   $$
   V^{2g-1}=\ker\left(H^{2g}(\bV^{<(2g-1)}))=
   \AA_c^{2g} \too H^{n+1}(\cN_X)\right) =
   ({\rm kernel} (F_c))^{2g}=\la q_g^1 \ra \cong \G_0
   $$
with $d\colon V^{2g-1} \isom \la q_g^1 \ra \subset \bV$.

For $n=2g$, we have
$$
H^{2g+1}(\bV^{\leq
(2g-1)})\,=\,H^{2g+1}(\bV^{<(2g-1)})\,=\,\AA_c^{2g+1}\, ,
$$
since $(\bV^{\leq
(2g-1)})^{2g+1}=\AA_c^{2g+1}\oplus (V^{2g-1}\cdot V^2)$ and the
nonzero elements in $V^{2g-1}\cdot V^2$ are not closed. So
   $$
   V^{2g}=({\rm kernel} (F_c))^{2g+1}=q_{g-1}^1
\cdot W_c \cong \G_{{\mathbf{e}}_1}\, ,
   $$
and $d\colon V^{2g} \isom q_{g-1}^1 \cdot V^3 \subset \bV$.

For $n=2g+1$, we have
  \begin{eqnarray*}
  (\bV^{\leq 2g})^{2g+1}&=&\AA_c^{2g+1}\oplus (V^{2g-1}\cdot V^2) \, , \\
  (\bV^{\leq 2g})^{2g+2}&=&\AA_c^{2g+2}\oplus
  (V^{2g-1}\cdot V^3) \oplus (V^{2g}\cdot V^2)\, ,
  \end{eqnarray*}
with $d\colon V^{2g-1}\cdot V^3 \isom q_g^1\cdot V^3$ and $d\colon
V^{2g}\cdot V^2\isom \alpha \ q_{g-1}^1 \cdot V^3$. But $\alpha \
q_{g-1}^1$ and $q_g^1$ are linearly independent, so we have
  $$
  H^{2g+2}(\bV^{\leq 2g})=\frac{Z^{2g+2}(\bV^{\leq 2g})}{B^{2g+2}(\bV^{\leq
  2g})}=\frac{\AA_c^{2g+2}}{\la \a \ q_g^1\ra} \, .
  $$
This gives
   $$
  V^{2g+1}=\ker \left( \AA_c^{2g+2}/\la \a q_g^1\ra \too
  H^{2g+2}(\cN_X)\right)
  = \la q_g^2 \ra \oplus q_{g-2}^1
  \cdot \bw^2_0 W_c \cong \G_0 \oplus \G_{{\mathbf{e}}_2}\, ,
   $$
which follows easily using Lemma \ref{lem:kerF}. Note that
the differential $d$ maps
the summand $\G_{{\mathbf{e}}_2}$ to $q_{g-2}^1 \cdot \bw^2_0
V^3$.

We now proceed to prove the second part of the theorem.

\textbf{Proof of (i):}
Let us prove by induction on $k\geq 1$ that there exists a
sub--representation $U_{k}\,\subset\, V^{2g+2k-2}$ with
$U_{k}\cong \G_{k \, {\mathbf{e}}_1}$ such that $d(U_{k})\subset
U_{k-1}\cdot V^3 \subset V^{2g+2k-4} \cdot V^3$, where $U_0:=\la
q_{g-1}^1 \ra\subset V^{2g-2}$.

If $k\,=\,1$, then
 $$
 V^n\,=\, V^{2g}\,=\,\G_{{\mathbf{e}}_1}
 $$
with $d\colon V^{2g} \,\too\,q_{g-1}^1 \cdot V^3 \subset \bV$.

Now assume that there exists a sub--representation
$U_{k-1}\,\subset\, V^{2g+2k-4}$ with $U_{k-1}\cong \G_{(k-1)
{\mathbf{e}}_1}$ such that $d(U_{k-1})\subset U_{k-2}\cdot V^3
\subset V^{2g+2k-6} \cdot V^3$. Then we have
$$
d(U_{k-1} \cdot V^3)\,\subset\, U_{k-2}\cdot V^3\cdot V^3\, .
$$
On one hand,
   $$
   \begin{array}{l}
   U_{k-1} \cdot V^3 =U_{k-1} \ox  V^3 \subset
   V^{2g+2k-4}\ox V^3\ , \\
   U_{k-1} \ox V^3 \cong \G_{(k-1){\mathbf{e}}_1} \ox
   \G_{{\mathbf{e}}_1}\, ;
   \end{array}
   $$
so there exists $\widetilde{U}_k\subset U_{k-1} \cdot V^3$ with
$\widetilde{U}_k \cong \G_{k\, {\mathbf{e}}_1}$ such that
$\widetilde{U}_k\subset U_{k-1}\cdot V^3$. On the other hand,
  $$
  U_{k-2}\cdot V^3\cdot V^3 \subset U_{k-2}\ox \bw^2 V^3 \cong
  \G_{(k-2){\mathbf{e}}_1} \ox \G_{{\mathbf{e}}_2}
  $$
does not contain $\G_{k \, {\mathbf{e}}_1}$. So
$d(\widetilde{U}_k)=0$, or in other words, $\widetilde{U}_k\subset
Z^{2g+2k-1}(\bV)$.

By Remark \ref{rem:extra}, we have $\rho(V^n)=0$ for all $n\geq
5$, where $\rho\colon (\bV,d)\too (H^*(\cN_X),0)$ is the minimal
model.  As ${U}_{k-1} \subset V^{2g+2k-4}$, we have
$\rho(U_{k-1})=0$. Hence $\rho(\widetilde{U}_k)=0$, or in other
words, $\widetilde{U}_k\subset B^{2g+2k-1}(\bV)$. This is only
possible if there exists $U_k\subset V^{2g+2k-2}$ with $U_k\cong
\G_{k\, {\mathbf{e}}_1}$ and $d\colon U_k\isom \widetilde{U}_k
\subset U_{k-1}\cdot V^3$. Therefore, the proof of
statement (i) is complete by induction.

\textbf{Proof of (ii):}
We will show using induction on $k\geq 1$ that there exists a
sub--representation
$$
U_{k}\,\subset\, V^{2g+2k+1}
$$
with
$U_{k}\cong \G_{k \, {\mathbf{e}}_1+ {\mathbf{e}}_2}$ such that
$d(U_{k})\subset U_{k-1}\cdot V^3 \subset V^{2g+2k-1} \cdot V^3$,
where $U_0\subset V^{2g+1}$ is the sub--representation of
$V^{2g+1}=\Gamma_0\oplus \Gamma_{\mathbf{e}_2}$ isomorphic to
$\Gamma_{\mathbf{e}_2}$.

For $k\,=\,1$, note that $d\colon U_0 \to q_{g-2}^1 \cdot \bw^2_0
V^3$. So $d\colon U_0\cdot V^3 \to q_{g-2}^1 \cdot \bw^2_0
V^3\cdot V^3$, where
 $$
 \bw^2_0 V^3\cdot V^3 \cong \bw^3_0 V^3 \oplus \gamma\cdot V^3 \cong
 \bw^3_0 W_c \oplus W_c =\Gamma_{\mathbf{e}_3}\oplus
 \Gamma_{\mathbf{e}_1}.
 $$
As $U_0\cong \Gamma_{\mathbf{e}_2}$, we conclude that
 $$
 U_0\cdot V^3\,\cong\, U_0\ox V^3\,\cong\, \Gamma_{\mathbf{e}_2}\ox
 \Gamma_{\mathbf{e}_1}
 $$
contains a sub--representation $\widetilde{U}_1\subset U_0\cdot
V^3$ with $\widetilde{U}_1\cong
\Gamma_{\mathbf{e}_1+\mathbf{e}_2}$ and $d(\widetilde{U}_1)=0$.
Working as in the proof of (i),
this yields that there exists $U_1\subset
V^{2g+3}$ with $U_1\cong \G_{ {\mathbf{e}}_1+ {\mathbf{e}}_2}$
and $d(U_1)=\widetilde{U}_1\subset U_0\cdot V^3$.

Now assume that $k\geq 2$ and that there exists a
sub--representation $U_{k-1}\,\subset\, V^{2g+2k-1}$ with
$U_{k-1}\cong \G_{(k-1) {\mathbf{e}}_1+{\mathbf{e}}_2}$ such that
$d(U_{k-1})\subset U_{k-2}\cdot V^3 \subset V^{2g+2k-3} \cdot
V^3$. Then we have $d(U_{k-1} \cdot V^3) \subset U_{k-2}\cdot
V^3\cdot V^3$. On one hand,
   $$
   U_{k-1} \cdot V^3 \,=\, U_{k-1} \ox  V^3
    \,\cong\, \G_{(k-1){\mathbf{e}}_1+\mathbf{e}_2} \ox
   \G_{{\mathbf{e}}_1}\, ;
   $$
so there exists $\widetilde{U}_k\subset U_{k-1} \cdot V^3$ with
$\widetilde{U}_k \cong \G_{k\, {\mathbf{e}}_1+\mathbf{e}_2}$ such
that $\widetilde{U}_k\subset U_{k-1}\cdot V^3$. On the other hand,
  $$
  U_{k-2}\cdot V^3\cdot V^3\,\subset\, U_{k-2}\ox \bw^2 V^3\,\cong\,
  \G_{(k-2){\mathbf{e}}_1+\mathbf{e}_2} \ox \G_{{\mathbf{e}}_2}
  $$
does not contain $\G_{k \, {\mathbf{e}}_1+\mathbf{e}_2}$. Therefore,
$d(\widetilde{U}_k)=0$. Thus there exists $U_k\subset V^{2g+2k-2}$
with $U_k\cong \G_{k\, {\mathbf{e}}_1+\mathbf{e}_2}$ and $d\colon
U_k\isom \widetilde{U}_k \subset U_{k-1}\cdot V^3$. This completes
the proof of the theorem.
\end{proof}

{}From Theorem \ref{thm:g>2} and Theorem
\ref{thm:leading} it follows that for each $g\, \geq\, 2$,
the rational homotopy group
$\pi_n(\cN_X)\ox_{\mathbb Z} \CC$
is nonzero for infinitely many $n$. As noted in the introduction,
this means that the moduli space $\cN_X$ is rationally hyperbolic
for all $g\, \geq\, 2$.
Therefore,
 \[
 f(k)\,=\, \sum_{i=1}^{\dim_{\mathbb R}\cN_X-1}
 \dim \pi_{k+i}(\cN_X)\ox_{\mathbb Z}
 {\mathbb Q}
 \]
grows faster than any polynomial in $k$.

\begin{remark}\label{rem:invar}
Let $X$ be a smooth irreducible projective complex curve of genus
$g\,\geq \, 2$. Whereas the minimal model of $(H^*(\cN_X,\CC),0)$
has infinitely many $n\in \NN$ for which $V^n\neq 0$, the minimal
model of the algebra $(H^*_I(\cN_X,\CC),0)$ has a very different
behavior. Actually, from eqn. \eqref{eqn:invar} we find that the
minimal model of $(H^*_I(\cN_X,\CC),0)$ is
  $$
  (\bbw (\alpha,\beta,\gamma, f_1,f_2,f_3),d), \quad df_1=q_g^1,\
  df_2=q_g^2,\ df_3=q_g^3),
  $$
where $\deg(\a)=2$, $\deg(\b)=4$, $\deg(\g)=6$, $\deg(f_1)=2g-1$,
$\deg(f_2)=2g+1$ and $\deg(f_3)=2g+3$.
\end{remark}

\noindent \textbf{Acknowledgements.}\, We are grateful to Aniceto
Murillo for useful comments. The first--named author wishes to
thank the Harish--Chandra Research Institute for its hospitality.
The second--named author is partially supported through grant MCyT
(Spain) MTM2004-07090-C03-01.

\end{document}